\documentclass[11pt]{article}
\usepackage[top=1in,
    bottom=1in,
    left=0.6in,
    right=0.6in
]{geometry}

\usepackage[T1]{fontenc}
\usepackage{microtype}
\usepackage{enumerate}
\usepackage{setspace}
\usepackage{listings}
\usepackage{amsmath,amsfonts,amssymb,graphics,amsthm,comment, bm, svg}
\usepackage{hyperref}
\hypersetup{
    colorlinks=true,
    linkcolor=blue,
    citecolor=blue,
    urlcolor=blue,
    pdfborder={0 0 0}
} 
\usepackage{dsfont}
\usepackage{tocloft}

\usepackage[framemethod=TikZ]{mdframed}
\usepackage{lipsum}
\mdfdefinestyle{MyFrame}{%
    linecolor=blue,
    outerlinewidth=2pt,
    roundcorner=20pt,
    innertopmargin=\baselineskip,
    innerbottommargin=\baselineskip,
    innerrightmargin=20pt,
    innerleftmargin=20pt,
    backgroundcolor=gray!50!white}
\usepackage{shadethm} 

\usepackage{graphicx}
\usepackage{float}
\usepackage{xcolor}
\usepackage{wrapfig} 

\usepackage[
    font={small,sf},
    labelfont={sf,bf},
    margin=1cm
]{caption}

\usepackage{longtable}
\usepackage{array}

\usepackage{cleveref}

\usepackage[backend=biber,style=alphabetic]{biblatex}
  \Crefname{theorem}{Theorem}{Theorems}
  \Crefname{thm}{Theorem}{Theorems}
  \Crefname{lemma}{Lemma}{Lemmas}
  
  \Crefname{lem}{Lemma}{Lemmas}
  \Crefname{remark}{Remark}{Remarks}
  \Crefname{prop}{Proposition}{Propositions}
  \Crefname{proposition}{Proposition}{Propositions}
  \Crefname{problem}{Problem}{Problems}
\Crefname{notation}{Notation}{Notations}
\Crefname{shadethm}{Theorem}{Theorems}
\Crefname{claim}{Claim}{Claims}
  \Crefname{defn}{Definition}{Definitions}
  \Crefname{corollary}{Corollary}{Corollaries}
  \Crefname{shadedefn}{Definition}{Definitions}
  \Crefname{section}{Section}{Sections}
  \Crefname{figure}{Figure}{Figures}
  \Crefname{shadeprop}{Proposition}{Propositions}
  \Crefname{exercise}{Exercise}{Exercises}
    \Crefname{assumption}{Assumption}{Assumptions}

\newtheorem{thm}{Theorem}[section]
\newshadetheorem{shadethm}{Theorem}[section]
\definecolor{shadethmcolor}{HTML}{BEFF33}
\definecolor{shaderulecolor}{HTML}{45CFFF}
\newshadetheorem{shadedefn}{Definition}[section]
\definecolor{shadethmcolor}{HTML}{CBFAFA}
\definecolor{shaderulecolor}{HTML}{CBFAFA}
\newshadetheorem{shadeprop}{Proposition}[section]
\definecolor{shadethmcolor}{HTML}{CAFAF2 }
\definecolor{shaderulecolor}{HTML}{CBFAFA}
\newtheorem{lemma}[thm]{Lemma}

\newtheorem{corollary}[thm]{Corollary}

\newtheorem{defn}[thm]{Definition}

\numberwithin{equation}{section}

\theoremstyle{definition}
\newtheorem{remark}[thm]{Remark}

\def \min {\text{min}}
\makeatletter
\def\munderbar#1{\underline{\sbox\tw@{$#1$}\dp\tw@\z@\box\tw@}}
\makeatother

\def\cP{\mathcal{P}}

\def\cN{\mathcal{N}}
\def\cM{\mathcal{M}}

\def\cE{\mathcal{E}}

\def\cC{\mathcal{C}}
\def\cB{\mathcal{B}}
\def\cA{\mathcal{A}}

\usepackage{wasysym}

\def\P{\mathbb{P}}
\def\E{\mathbb{E}}

\def\N{\mathbb{N}}

\def  \p- {p\textunderscore}

\def\eps{\varepsilon}

\def \Kn {\vec K_n}

\DeclareMathOperator{\var}{Var}

\DeclareMathOperator{\cbd}{\cC^{\partial}}
\DeclareMathOperator{\cmax}{\cC^{\sf{max}}}

\DeclareMathOperator{\Var}{Var}
\DeclareMathOperator{\Cmax}{{\sf{Cycle}}_{\max}}

\usepackage{extarrows}

\title{Central Limit Theorem for the Weight of the Minimum Spanning Arborescence on the Complete Directed Graph}
\author{Swarnadeep Bagchi \thanks{University of Victoria, Victoria, Canada. Research was supported in part by NSERC. Email: swarnadeepbagchi@uvic.ca} \and Pablo de la Fuente\thanks{University of Victoria, Victoria, Canada. Research was supported in part by NSERC and SECIHTI. Email: pablodelafuente@uvic.ca} \and Kesav Krishnan \thanks{Indian Institute of Technology Madras, India. Research was supported in part by ANRF-ECRG. Email: kesav@iitm.ac.in}}
\date{ }

\begin{document}
\maketitle
\abstract{The minimum spanning arborescence (MSA) with a fixed boundary, defined on weighted directed graphs, is the directed analogue of the minimum spanning tree (MST). The behaviour of the MSA has been of recent interest and is known to be quite different from that of its undirected counterpart. On the complete directed graph, and for edge weights given by i.i.d.\ random variables with a broad class of distributions, we prove a Central Limit Theorem for the weight of the MSA with appropriate scaling. While the scaling for the CLT is similar to MST case explored by Janson, the centering required ($1$ as opposed to $\zeta(3)$ for Uniform$[0,1]$ weights) and the limiting variance of the Gaussian ($1$ as opposed to $6\zeta(4)-4\zeta(3)$ for Uniform$[0,1]$ weights)  are notably different. }

\tableofcontents

\section{Introduction}
Let $\vec G=(V,\vec E)$ be a graph where each edge comes with a specific direction and there are no self-loops. More precisely, $\vec E$ is a subset of the following set of ordered pairs
\[
\{(u,v): u,v\in V, u\neq v\}.
\]
If $\vec e\in\vec E$, such that $\vec e=(u,v)$, we say that $u$ is the \emph{tail} of $\vec e$ and $v$ is the \emph{head} of $\vec e$ and denote them by $\vec e_{-}$ and $\vec e_{+}$, respectively. We say $\vec e$ is an outgoing edge from $u$ and an incoming edge to $v$. A \emph{directed path} is a sequence of distinct edges $\vec e_1, \vec e_2,\cdots,\vec e_n$ such that ${\vec e_{i_+}}={\vec e_{(i+1)_{-}}}$, for all $1\le i < n$. If there is a directed path $\vec e_1, \vec e_2,\cdots,\vec e_n$ such that ${\vec e_{n_{+}}}={\vec e}_{1_{-}}$, we say that this path creates a \emph{directed cycle}. Sometimes we denote a directed cycle by the set of heads and tails of the edges in it. If there is a directed path from a vertex $u$ to a vertex $v$, we say that $u$ is in the \emph{past} of $v$, and denote it as $u\in{\sf{Past}}_v$. Note that our definition admits paths with no edges, so $v\in{\sf{Past}}_v$, for all $v\in V$.\\

\begin{figure}[H]
\centering
\includegraphics[width=12cm]{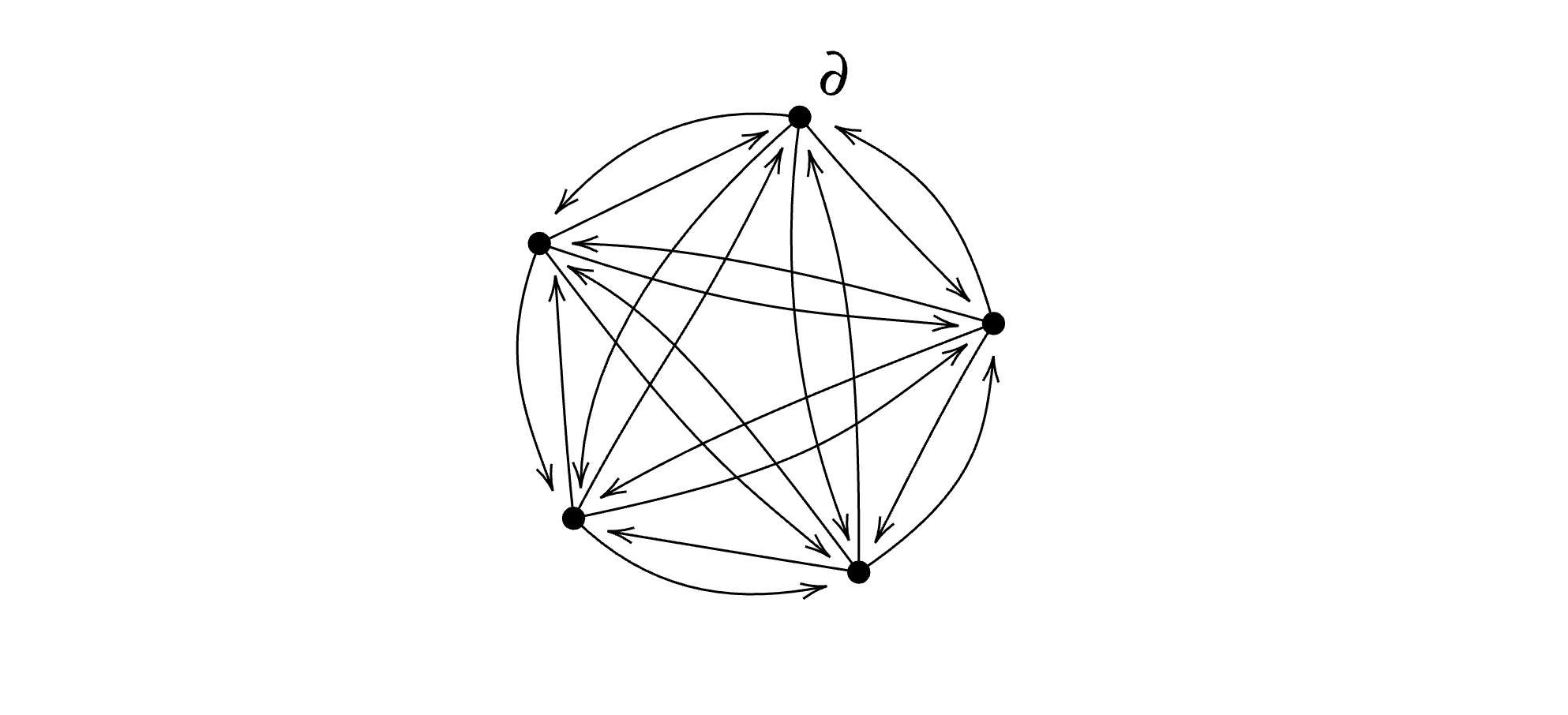}
\caption{The complete directed graph on $5$ vertices with distinguished vertex $\partial$.}
\label{fig:directedkn}
\end{figure}

The setting of this article is on the complete directed graph $\vec{K}_{n}=(V_{n},\vec{E}_{n})$, where $V_{n}$ is a set of $n$ vertices, and $\vec{E}_{n}$ is the set of all possible directed edges. More precisely, the edge set of the complete directed graph is given by  
\[
\vec{E}_{n}=\{\vec{e}=(u,v):u,v\in V_n,\ u\neq v\}.
\]
See Figure \ref{fig:directedkn} for an example.

\begin{defn} \label{def:arborescence}
An \textbf{arborescence}, denoted $(\vec T,\partial)$ is a finite directed graph $\vec T=(V(\vec T), \vec E(\vec T))$ with boundary $\partial\in V(\vec T)$, which has the following properties:
\begin{enumerate}
    \item $\partial$ has no outgoing edge,
    \item all vertices $v\in V(\vec T)\setminus\{\partial\}$ have exactly one outgoing edge,
    \item there are no directed cycles and
    \item the graph is connected in the undirected sense.    
\end{enumerate}
\end{defn}

\begin{defn}
    A \textbf{spanning arborescence} $(\vec T,\partial)$ of a finite directed graph $(\vec G, \partial)$, where $\partial$ is a distinguished vertex, is an arborescence with boundary $\partial$ such that $V(\vec{T})=V(\vec G)$ and $E(\vec{T})\subset E(\vec{G})$. 
\end{defn}

\begin{figure}[H]
    \centering
    \includegraphics[width=12cm]{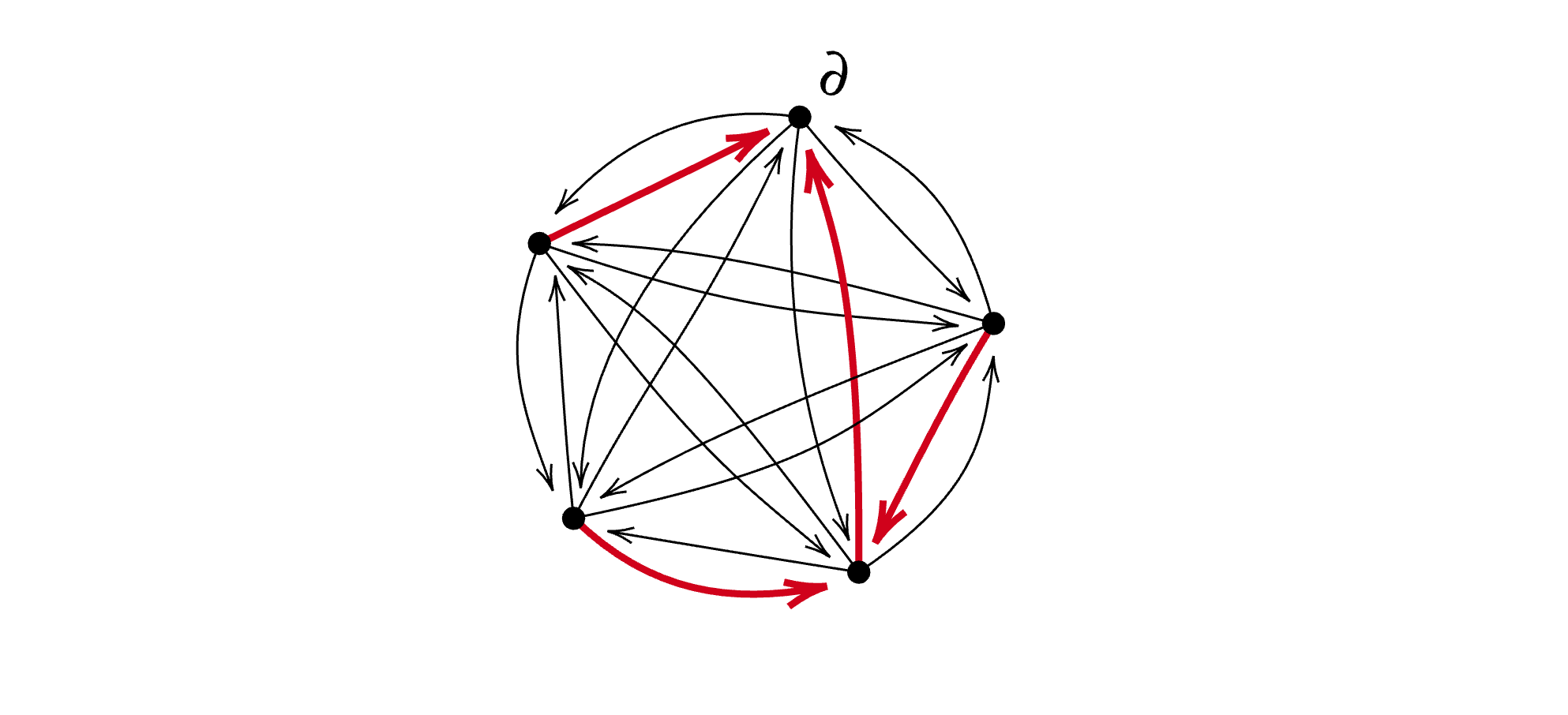}
    \caption{The red edges form a spanning arborescence of $(\vec K_5,\partial)$.}
    \label{fig:arborescence}
\end{figure}

See Figure \ref{fig:arborescence} for an example.

\begin{remark}
    For a spanning arborescence, given conditions 1 and 2 in Definition \ref{def:arborescence}, points 3 and 4 are equivalent to saying that the underlying undirected graph is a spanning tree, and then all edges are directed towards $\partial$. 
\end{remark}
\begin{remark}
In an arborescence, from each vertex there is a unique directed path to $\partial$. 
\end{remark}
The article concerns the behaviour of the spanning arborescence of $(\vec{K}_{n}, \partial)$ that has minimum weight when the edges are assigned random weights. This means that each of the directed edges of $\Kn$ is assigned an independent non-negative random variable with density $p(x)$, which satisfies the following properties: 
\begin{itemize}
\item $p(0)>0$,
\item $p$ is continuous at $0$,
\item $\displaystyle\int_0^{\infty}xp(x)dx<\infty$.
\end{itemize}
We denote this family of random variables by $\{Y_{\vec{e}}\}_{\vec{e}\in \vec{E}_{n}}$
and define the weight of a spanning arborescence $(\vec T,\partial)$ of $(\Kn,\partial)$ to be 
\[
W(\vec{T})=\sum_{\vec{e}\in \vec E(\vec T)}Y_{\vec{e}}.
\]
The minimum spanning arborescence (MSA) on $(\vec{K}_{n},\partial)$ is then defined to be 
\[
(\vec{T}^{*}_{n},\partial):=\text{argmin} \{W(\vec{T}):(\vec{T},\partial) \text{ is a spanning arborescence on }(\Kn,\partial)\}.
\]
We will denote by $W_{n}^{*}$ the weight of the minimum spanning arborescence. Note that the minimum spanning arborescence is almost surely unique as the weights on the edges are continuous random variables. The main result of this article is the following: 
\begin{thm}\label{thm:main}

Let i.i.d.\ weights $\{Y_{\vec e}\}_{\vec e\in\vec E_n}$ be assigned to the edges of $(\Kn, \partial)$ whose common density $p$ satisfies the previous hypotheses. Let $X^{(n)}_v:=\min_{\vec e_{-}=v}Y_{\vec e}$. Then as $n\to\infty$,
\[
\sqrt n\left(W_n^*-n\E\left(X^{(n)}_v\right)\right)\xrightarrow[]{d}\cN\left(0,\frac{1}{p(0)^2}\right).
\]
\end{thm}
Since $\{Y_{\vec e}\}_{\vec e\in\vec E_n}$ are i.i.d.,\ $\{X^{(n)}_v\}_{v\in V_n\setminus\{\partial\}}$ are also i.i.d. For brevity we remove the suffix and denote them as $X^{(n)}$, for all $v\in V_n\setminus\{\partial\}$. The limiting variance is $\frac{1}{p(0)^2}$ since $\lim_{n\to \infty}n^2\var \left( X^{(n)} \right) = \frac{1}{p(0)^2}$, we can derive it as an immediate corollary of Lemma \ref{lem:genweight}.

\begin{remark}
Let $Z_n=\sqrt n\left(W_n^*-n\E\left(X^{(n)}_v\right)\right)$ and $a_n=n\E\left(X^{(n)}_v\right)$. By Theorem \ref{thm:main}, we know that $Z_n$ is tight. Since $n^{-1/2}\to 0$, $W_n^*-a_n=Z_n/\sqrt n\to 0$ in probability. If the edge weights are i.i.d. Exponential$(1)$ or Uniform$[0,1]$, then $a_n=1$. Thus $W_n^*\to 1$ in probability.
\end{remark}
Our article follows a great deal of recent interest in understanding the behaviour of the MSA in various settings, including the complete directed graph. The issue of almost sure convergence to an infinite volume limit was explored in \cite{RS} by Ray and Sen. Inspired by an algorithmic construction of the MSA, they describe a new stochastic process termed the ``loop-contracting random walk'' for the purposes of constructing the arborescence. This is in analogy with the loop-erased random walk and its relationship with weighted spanning trees, formalized by Wilson's algorithm. One of the major results of \cite{RS} is that the infinite volume limit exists almost surely if the loop-contracting random walk is almost surely transient. Verifying the transience is tricky since the process is not Markovian and it is unclear if it admits a zero-one law.

The MSA on the complete directed graph has been studied in \cite{Frieze1985}, \cite{FriezeTkocz2022} and \cite{BR}. In the undirected case, a law of large numbers for the weight of the minimum spanning tree was proved by Frieze \cite{Frieze1985}, remarkably showing that for a wide class of weights the limit is $\zeta(3)$. The method of \cite{Frieze1985} hinges on assuming that the density of the weights satisfies $p(0)>0$. For the MSA, the analogous law of large numbers for uniform $[0,1]$ weights was accomplished by Frieze and Tkocz in \cite{FriezeTkocz2022}. The limit, in contrast to the undirected version's $\zeta(3)$, is $1$.  The central limit theorem for the weight of the MST on the complete graph was proved by Janson \cite{janson1995minimal} with i.i.d. Uniform$[0,1]$ edge weights. The variance of the limiting Gaussian was first identified as a highly complicated sum, later simplified by W\"astlund to $6\zeta(4)-4\zeta(3)$ \cite{Wastlund2005}.

Hansen  established a CLT, as well as a law of large numbers (in the $L^{1}$ sense), for the weight of the MSA \textbf{without a fixed boundary} \cite{hansen1997limit}. Fixing a boundary vertex imposes additional connectivity constraints, making the repair steps (specifically {\sf Step 2} and  {\sf Step 3} detailed in  Section \ref{sec:proofmethod}) different. Furthermore, the CLT of \cite{hansen1997limit} can be recovered from our result for the class of weights considered in this article. We discuss this in Section \ref{sec:concreal}. The case with a designated boundary point is important, since there are known algorithms for constructing the rooted MSA, such as the Chu-Liu-Edmonds-Bock (CLEB) algorithm (see \cite{bock}, \cite{Edmonds}). Constructing the MSA without a designated boundary point is significantly more tedious, whether the rooted case can be bypassed altogether remains unclear. The CLEB algorithm is instrumental in the study of the infinite volume structure MSAs in \cite{RS} and \cite{BR}.  

Bagchi and Ray explored the limiting structure of the MSA as a graph and characterized the local weak limit as a Poisson$(1)$ Galton-Watson tree conditioned to survive directed towards an infinite path. To identify the local weak limit of MSA on the complete directed graph, the memoryless property of the exponential distribution plays very important role. This stands in sharp contrast to the local weak limit of the minimum spanning tree on the undirected complete graph. In both \cite{RS} and \cite{BR}, it is observed that the geometry of the MSA  depends on the distribution of the weights, significantly complicating its study. 

More generally, our article fits neatly into the theme of exploring central limit theorems for randomized combinatorial optimization problems. We direct the reader to \cite{Cao} and \cite{DT} for CLTs for certain other combinatorial optimization problems in the mean field setting. Note, \cite{DT} is in the positive temperature setting and can thus be thought of as CLTs for approximate versions of the relevant optimization problems. 
\subsection*{AI Declaration}
All ideas for the proof were conceived by us. Google Gemini Pro and ChatGPT were used for verification of details, finding gaps and proof reading. Gemini also suggested the integral representation used in the saddle point argument in the proof of Lemma \ref{lem:cyctail}. We subsequently verified all details by hand, and the actual writing has been done by us. 

\subsection*{Acknowledgement} The authors would like to thank Gourab Ray for suggesting the problem and for several subsequent useful discussions, concerning both the mathematics and presentation.

\section{Proof Method}\label{sec:proofmethod} In Sections \ref{sec:repair}, \ref{sec:structure} and \ref{sec:map}, we present detailed proofs for the case when weights have exponential distribution with mean $1$, mainly because the memoryless property of the exponential distribution allows several calculations in the proof to be carried out explicitly. We explain how the methods carry over to a more general family of distributions in Section \ref{sec:conc}.\\

Throughout the article, for positive $f$ and $g$, we are going to use the following standard asymptotic notation:
\begin{itemize}
    \item $f(n)=O(g(n))$ if there exist
constants $C>0$ and $n_0\in\mathbb{N}$ such that
$f(n)\le Cg(n)$
for all $n\ge n_0$.
\item $f(n)=\Theta(g(n))$ if there exist
constants $c,C>0$ and $n_0\in\mathbb{N}$ such that $cg(n)\le f(n)\le Cg(n)$,
for all $n\ge n_0$. If $f(n)/g(n)\to 1$, then we say that $f\sim g$. 
\item $f(n)=o(g(n))$ if
$\displaystyle\lim_{n\to\infty}\frac{f(n)}{g(n)}=0.$
\end{itemize}

\begin{figure}[H]
    \centering
    \includegraphics[width=8cm]{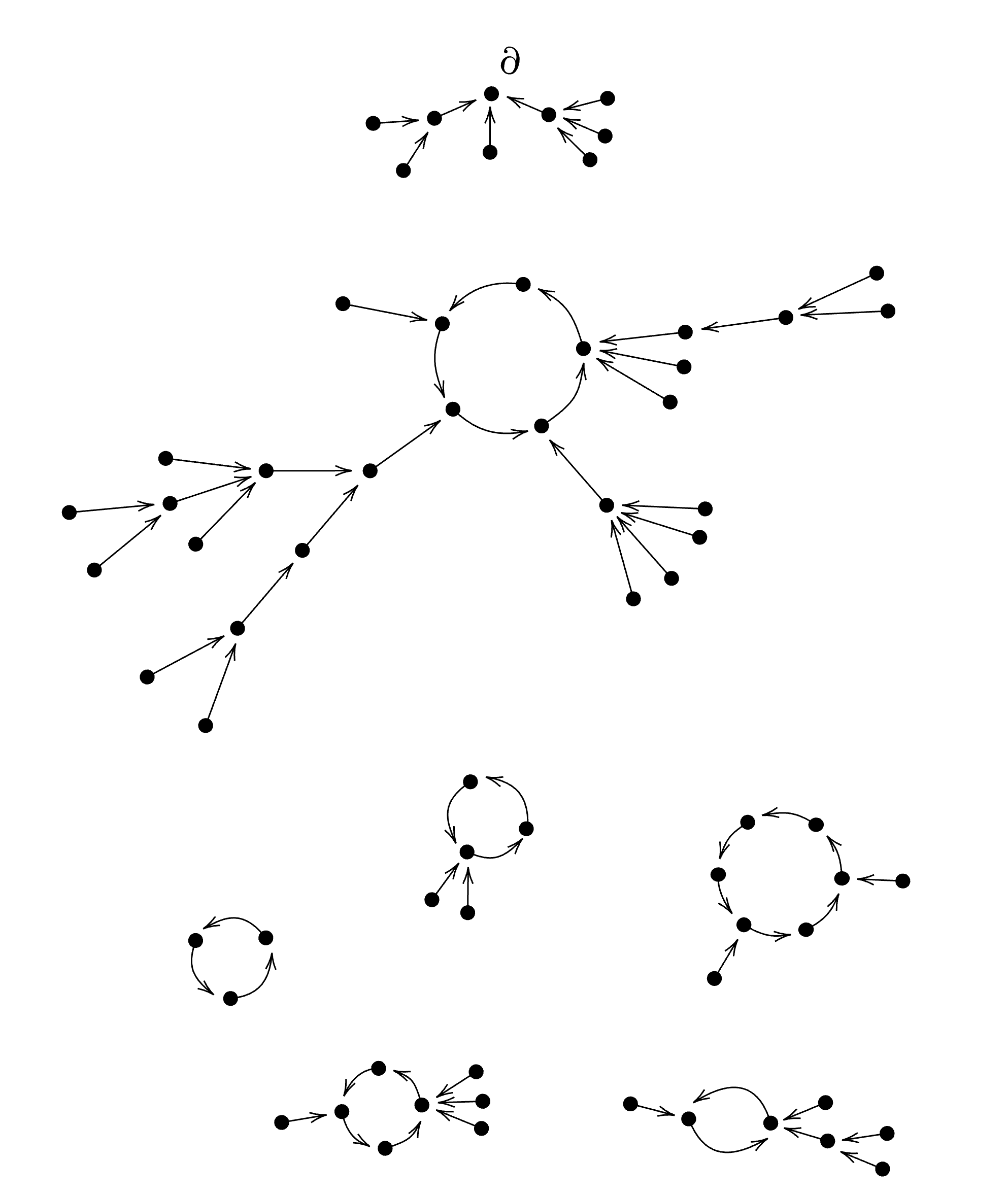}
    \caption{$\vec{K}_{n}^{1}$ obtained by exposing the minimum-weight outgoing edges from all the vertices, except $\partial$, in $\Kn$. Every component has one cycle, except $\cbd$ that is an arborescence with boundary $\partial$.}
    \label{fig:kn1}
\end{figure}

Recall that we have the complete directed graph $\Kn$ with a fixed vertex $\partial$ and every directed edge $\vec e\in\vec E_n$ independently has an Exponential$(1)$ weight $Y_{\vec e}$. The first step is to expose the minimum-weight outgoing edge from every vertex except $\partial$. We denote the configuration so obtained by $\vec{K}_{n}^{1}$, and its weight by $W_n$, that is $W_n=\sum_{\vec e\in\vec E(\Kn^1)}Y_{\vec e}$. See Figure \ref{fig:kn1} for a possible configuration.

Any two vertices that are connected by a sequence of edges, irrespective of direction, are said to be part of the same connected component. Note the following points about the connected components of $\Kn^1$, which are straightforward to prove: 
\begin{itemize}
    \item The connected component containing $\partial$ is an arborescence with boundary $\partial$. 
    \item Every other connected component consists of exactly one directed cycle with arborescences going into it.
\end{itemize}

We will comment more on the structural properties of this directed graph $\Kn^1$ in the next section. It suffices for now to note that the choice of exposed edges from each vertex, as well as their weights, are independent.

\begin{lemma}
    $W_n \rightarrow 1$ almost surely as $n\to\infty$.
\end{lemma}
\begin{proof}
    Let $X_v^{(n)}$ be the weight of the edge exposed from vertex $v$ when exposing the minimum weight edge. This means, 
    \[
    X_{v}^{(n)}:=\min_{\vec e_-=v} Y_{\vec{e}}.
    \]
    The minimum of $m$ independent exponential random variables is an exponential random variable whose mean is $m^{-1}$. We can then rescale these variables to get 
    \[
    X_v^{(n)}=\frac{Y_{v}^{(n)}}{n-1},
    \]
    where $\{Y_{v}^{(n)}\}_{v\in V_{n}}$ are independent exponential random variables with mean $1$.
    Since 
    \[
    W_n = \sum_{v\in V_n\setminus\{\partial\}} X_v^{(n)} = \frac{1}{n-1}\sum_{v\in V_n\setminus\{\partial\}}Y_v^{(n)},
    \]
    the almost sure limit $W_n \rightarrow \mathbb{E}(Y_v) = 1$ follows from the strong law of large numbers. \qedhere
 \end{proof}

 
 For the next layer of asymptotics, note that it is straightforward to extract a CLT for $W_{n}$, which we now state without proof:
 \begin{lemma}\label{lem:Wclt}
 \[
 \sqrt{n-1}\left( W_n - 1\right)\overset{d}{\rightarrow} \cN (0,1).
 \]
 \end{lemma}


Recall that $W_n^*$ is the weight of the MSA in $(\Kn, \partial)$. It turns out the law of large numbers for $W_{n}$ is an important first step towards verifying a law of large numbers for $W_{n}^{*}$. Thanks to that, the authors in \cite{FriezeTkocz2022} reduce their problem to proving that $W^*_n- W_n $ converges to $0$ almost surely. To start with, note that $W_n \leq W_n^*$, since every vertex (except $\partial$) has an outgoing edge and the weight of it is greater than or equal to the weight of the minimal edge. The technique is to alter $\vec{K}_{n}^{1}$  such that the resulting graph is a spanning arborescence of $(\Kn,\partial)$, which is guaranteed to have weight greater than or equal to $W_{n}^{*}$.  We refer to this procedure as ``repairing'', and denote its weight by $R_n$. Thus, $R_n$ is the sum of the weights of the edges added in by the modifications minus the weights of the edges removed. 

In \cite{FriezeTkocz2022} the law of large numbers for $W_{n}^{*}$ then follows from finding a repair procedure such that $R_{n}\to 0$ almost surely.  However, to show that the CLT for $W_{n}$ extends to $W_n^*$ we further need to show that $\sqrt{n-1} (W^*_{n}-W_{n})\to 0$ in probability, and then the extension follows from Slutsky's theorem. The problem is thus to find a repair procedure such that
\[
 \sqrt{n}R_{n}\to 0
\]
in probability.\\

The following terminology will be used for the remainder of the article. 
\begin{enumerate}
\item We denote by ${\sf Comp}$ the collection of connected components of $\vec{K}_{n}^{1}$
\item We denote by $\cmax$ the largest connected component of $\vec{K}_{n}^{1}$.
\item We denote by $\cbd$ the connected component containing $\partial$.
\item We denote by $\Cmax$ the directed cycle contained in $\cmax$.
\end{enumerate}

 We identify certain structural properties for the  configuration $K_{n}^{1}$, which if satisfied, admit a repair procedure yielding an arborescence of small enough weight. These properties are explicitly laid out in \eqref{eq:good}. To summarize:
\begin{enumerate}
\item The number of connected components of $\vec{K}_{n}^{1}$, $|{\sf Comp}|$ is bounded above by $C\log n$ for some constant $C<\infty$.
\item $|\Cmax|$ is bounded below and above by $|\cmax|^{\frac{1}{2}-\delta}$ and $|\cmax|^{\frac{1}{2}+\delta}$ respectively for some $\delta>0$. 
\item We have a large number of vertices $u\in \cmax$ such that ${\sf{Past}}_u$ is appropriately large.
\end{enumerate}
The event where these occur holds with high probability. This is precisely stated in Lemma \ref{lem:good} and a proof is presented in Section \ref{sec:map}. An important consequence for our repair procedure is that for $n$ sufficiently large (in a deterministic sense), on the good event, $\cbd\neq \cmax$. This is because the cycle size is large.   \\

\begin{figure}[H]
    \centering
    \includegraphics[width=7.7cm]{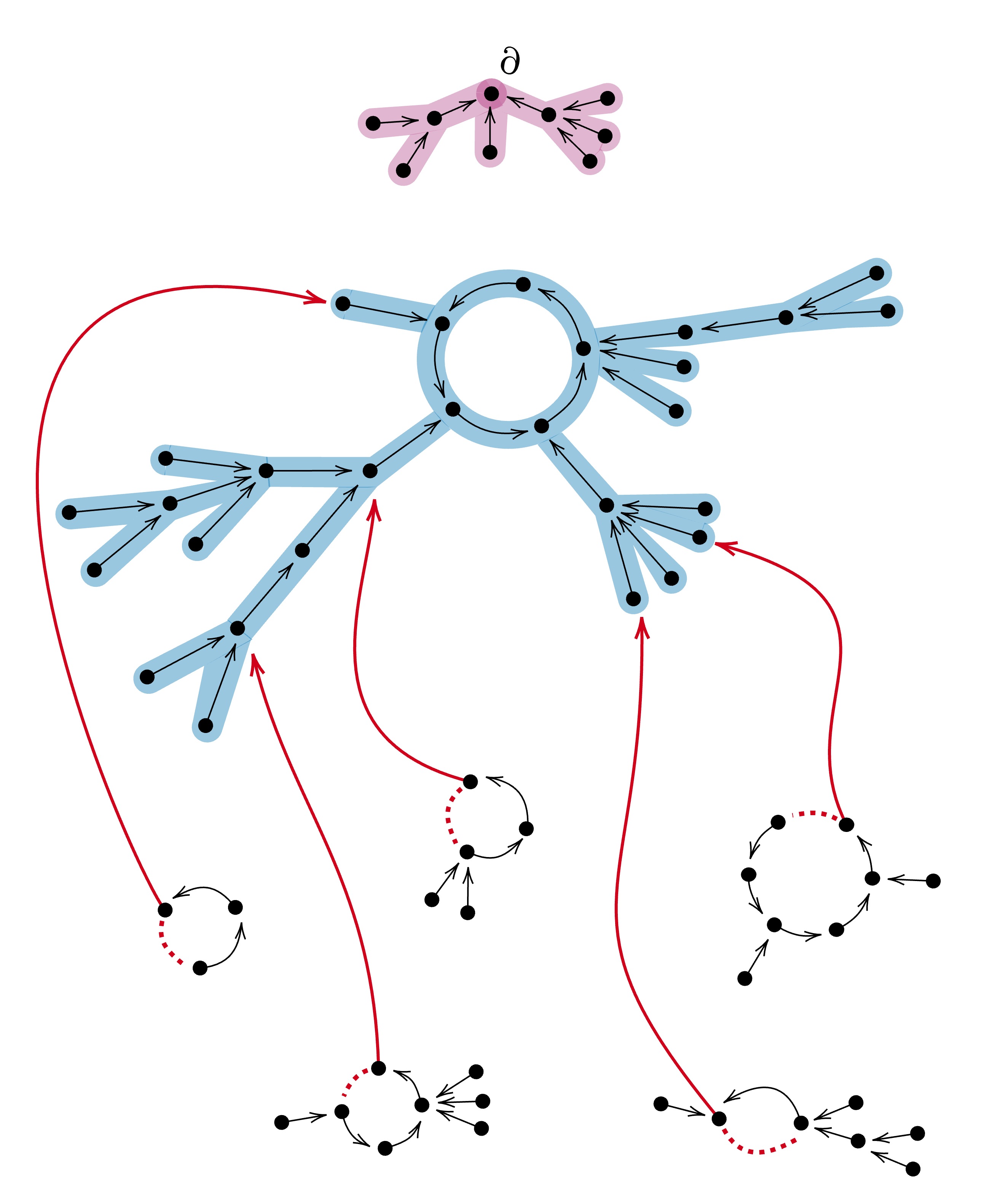}
    \caption{\textsf{Step 1} of the repair. From the vertices of the cycles of the components other than $\cbd$ (coloured pink), we pick the minimal outgoing edges (coloured red) to $\cmax$ (coloured blue), add these edges and delete the edges (dotted red) which are in the cycles.}
    \label{fig:step1}
\end{figure}

\begin{figure}[H]
    \centering

    \begin{minipage}[c]{0.5\textwidth}
        \vspace{-2em}
        \centering
        \includegraphics[width=\linewidth]{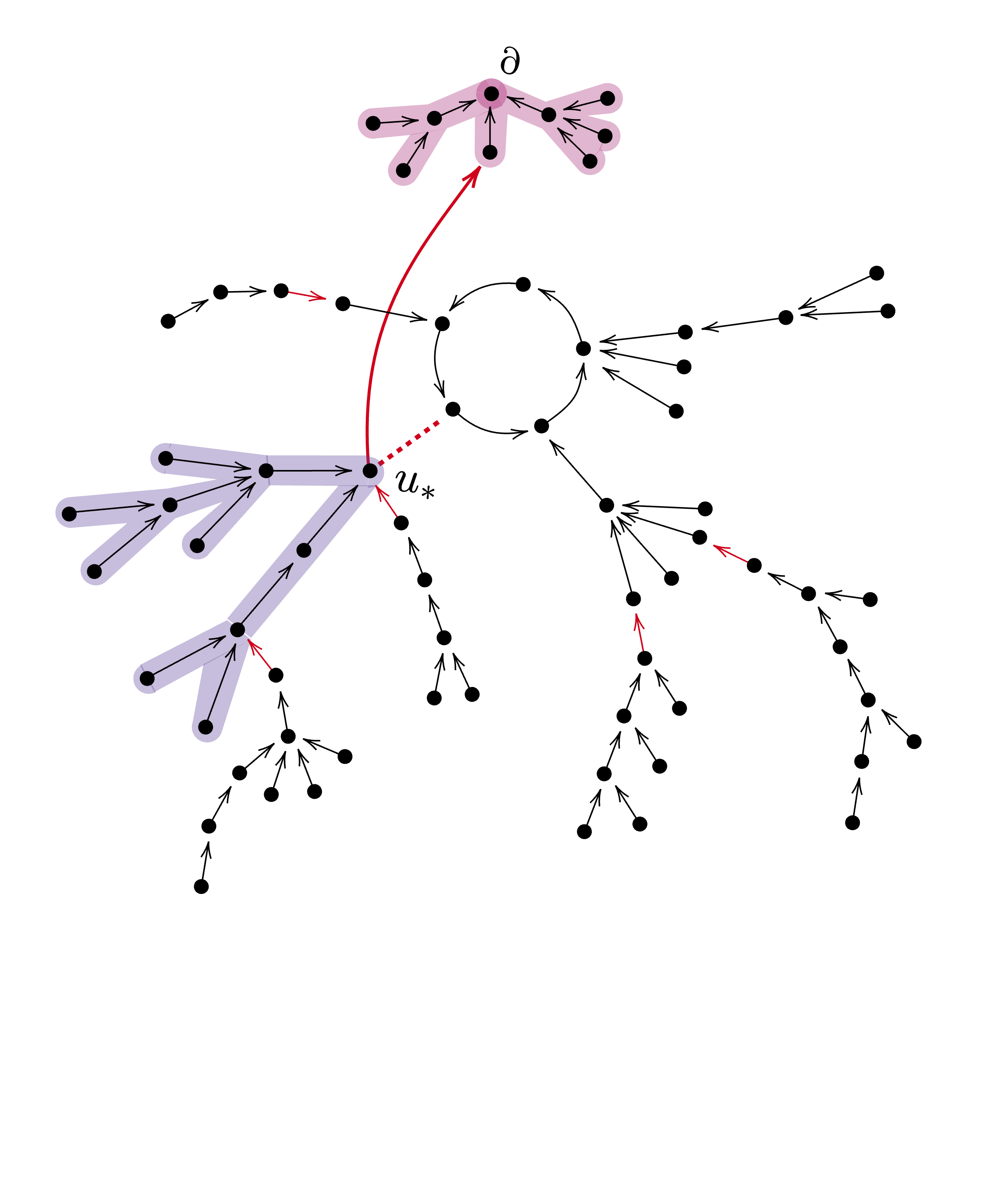}
    \end{minipage}%
    \hfill
    \begin{minipage}[c]{0.5\textwidth}
        \vspace{-2em}
        \centering
        \includegraphics[width=\linewidth]{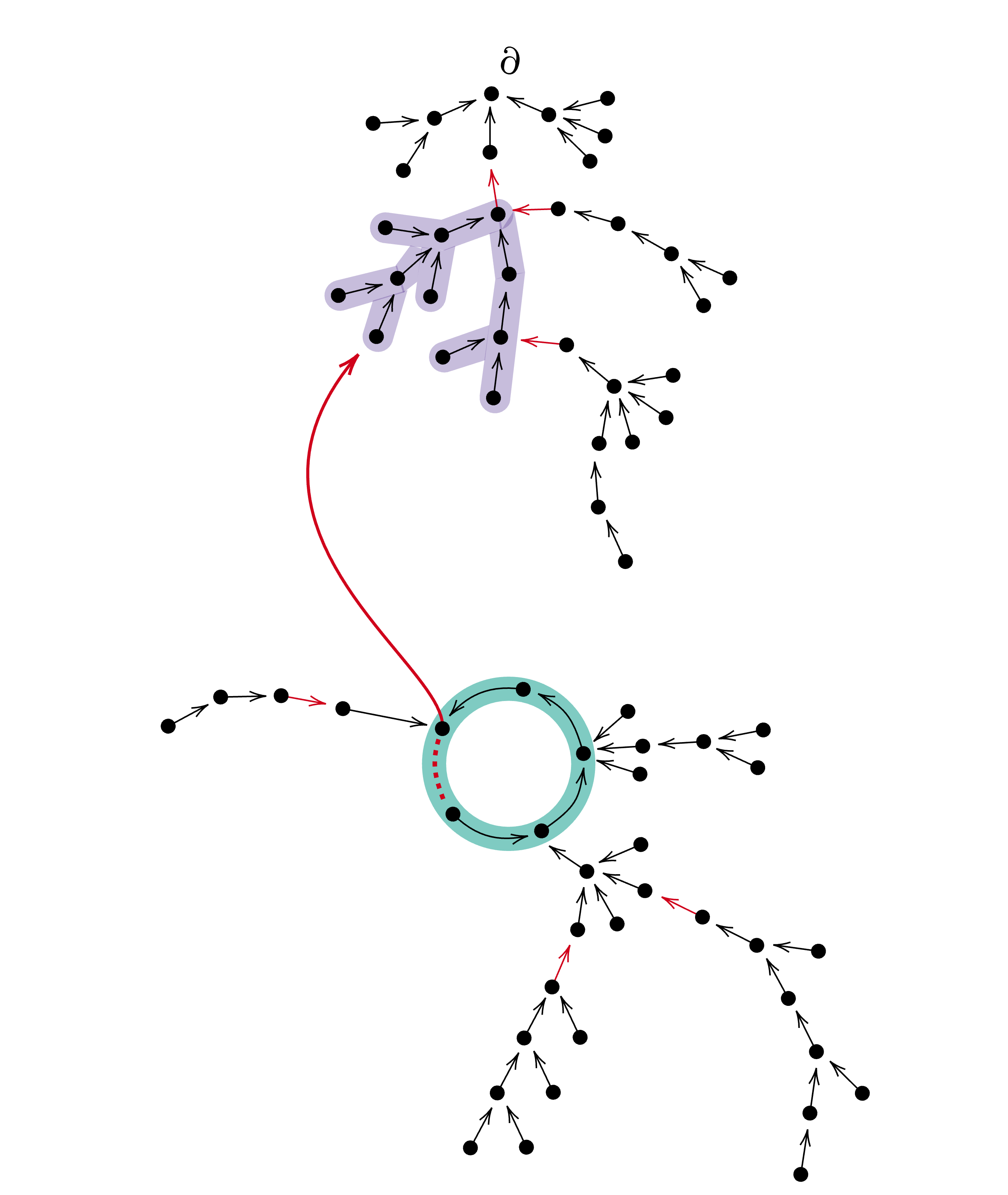}
    \end{minipage}

    \includegraphics[width=0.55\textwidth]{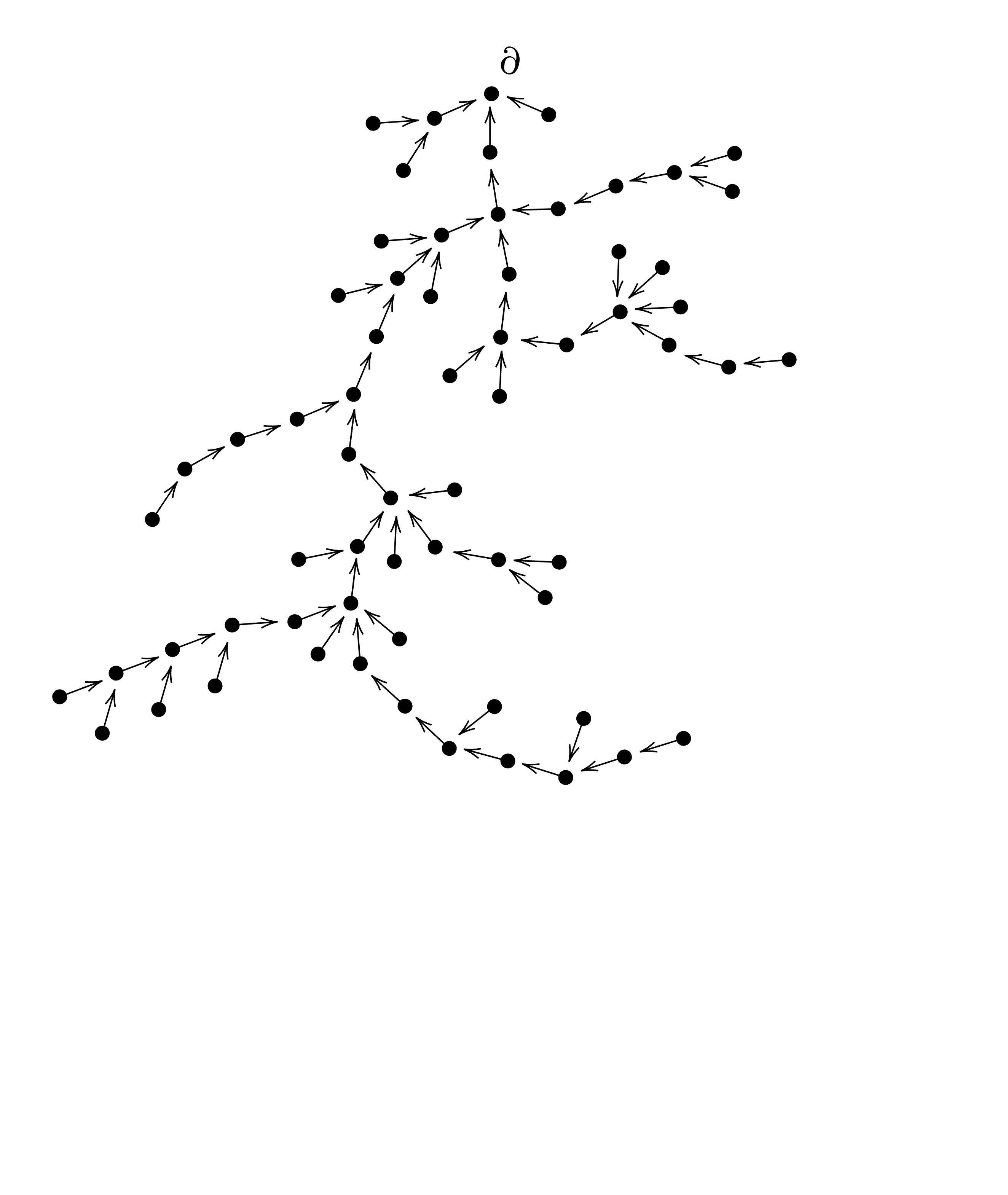}

    \vspace{-3cm}
    \caption{\textbf{Top left.} \textsf{Step 2} of the repair. $u_*$ is the vertex in $\cA\subset\cmax$ from which the minimal outgoing edge (coloured red) is hitting $\cbd$. We add this edge and delete the red dotted edge from $u_*$ towards the cycle. This step joins ${\sf{Past}}_{u_*}$ (coloured purple) to $\cbd$ (coloured pink).\\ \textbf{Top right.} \textsf{Step 3} of the repair. We pick the minimal outgoing edge (coloured red) from $\Cmax$ (coloured teal) to ${\sf{Past}}_{u_*}\cup\cbd$. Add that edge and deleted the red dotted edge in $\Cmax$.\\ \textbf{Bottom.} The spanning arborescence of $(\Kn,\partial)$ we obtain after finishing \textsf{Step 3} of the repair.}
    \label{fig:repair}

\end{figure}

Once we are on the good event, the repair procedure is as follows. For illustrations on the repair steps see Figures \ref{fig:step1} and \ref{fig:repair}.

\begin{itemize}
\item \textsf{\textbf{Step 1}}: Merge all components, except $\cbd$, with $\cmax$. This is done by adding an outgoing edge from the cycle in a given component $\cC$ to any vertex in $\cmax$, and from the vertex that now has two outgoing edges, delete the edge that leads back into the cycle. We choose the edge to be the one that minimizes the increase in weight of this repair, that is, the weight of the added edge minus the weight of the deleted edge. For brevity, we keep the name of $\cmax$ unchanged after this step.

The error accumulated in this step is bounded in Lemma \ref{lem:errbd1}. For the new edge from $\cC$ to $\cmax$, there are at least $|\cmax|$ many edges to choose from, the minimal edge over these having expected weight $O(|\cmax|^{-1})$. The lower bound on $|\cmax|$ follows from the upper bound on $|{\sf Comp}|$ via the pigeonhole principle. We then have to add the errors arising from all components apart from $\cbd$, yielding an overall expected error bound of $O(|\cmax|^{-1}|{\sf Comp}|)$.
\end{itemize}

At this step of the process, we have two disjoint arborescences. The obstacle now is merging $\cmax$ with $\cbd$ by connecting a vertex in $\Cmax$ to $\cbd$. The naive approach fails here. 
The best bounds we have are $|\cbd|\geq 1$ as $\cbd$ is going to be distributed as a Poisson$(1)$ Galton-Watson tree as $n\to\infty$, and $|\Cmax|\geq n^{\frac{1}{2}-\delta}$, for some small positive $\delta$, implying that all we can say about the minimal repairing edge is that it has weight $O(n^{-\frac{1}{2}+\delta})$. This does not tend to 0 when multiplied by $\sqrt{n}$. The issue is not just the quality of the bound, since in reality we can prove that $|\cbd|$ and $|\Cmax|$ are of order $1$ and $\sqrt{n}$ respectively, so there is no hope that this repair works even with an improved lower bound. The next step therefore relies intuitively on increasing the size of the boundary component, and thus reducing the weight of the minimal edge required for connection. It is for this purpose that the notion of the pasts of the vertices in $\cmax$ is required, as well as the fact that on the good event there are many of them with sufficiently large pasts.

\begin{itemize}
\item \textsf{\textbf{Step 2}}:  Denote by $\cA$ the set of all $u$ with sufficiently large pasts in $\cmax$. The precise notion can be found in Definition \ref{def:A}, with relevant quantifiers set up in \eqref{eq:epsdeletaprop}.

For each $u\in\cA$, let $E_{\partial}(u)$ be the set of all outgoing edges from $u$ to $\cbd$ in $\Kn$. We first select, for each $u\in\cA$, the edge in $E_{\partial}(u)$ that minimizes the increase in weight obtained by replacing the current outgoing edge from $u$. We then choose the vertex $u_*\in\cA$ for which this minimum increase is smallest, and perform the corresponding replacement.

The past of the chosen vertex, $\sf{Past}_{u_*}$, is an arborescence with boundary $u_*$. By deleting the original outgoing edge from $u_*$ and adding the selected edge in $E_{\partial}(u_*)$, we attach $\sf{Past}_{u_*}$ to $\cbd$, thereby increasing the size of $\cbd$. The increase in weight of this step has mean $O(|\cA|^{-1})$ and is handled precisely in Lemma \ref{lem:errbd2}.  

\item \textsf{\textbf{Step 3}}: We pick the outgoing edge from the vertices of $\Cmax$ to the vertices of ${\sf{Past}}_{u_*}\cup\cbd$ that minimizes the increase in weight of the repair. Suppose the tail of that edge is $u'$. We add that edge and delete the original outgoing edge emanating from $u'$ in $\Cmax$ in $\Kn^1$. This attaches $\cmax$ to $\cbd$, yielding a spanning arborescence of $\Kn$ with boundary $\partial$. The increase in weight from this step has mean $O(|\Cmax|^{-1}\cdot|{\sf{Past}}_{u_{*}}|^{-1})$, and is precisely bounded in Lemma \ref{lem:errbd3}. The key point is that the good event is specified such that we can ensure $|\Cmax||{\sf{Past}}_{u_*}| >>\sqrt{n}$, specifically $n^{\frac{1}{2}+\epsilon}$, for any $\epsilon>0$ as per \eqref{eq:epsdeletaprop}.  

\end{itemize}

Section \ref{sec:repair} concludes with the proof of Theorem \ref{thm:main}, modulo the proof of Lemma \ref{lem:good}. Additional structural ingredients required to prove Lemma \ref{lem:good} appear in Sections \ref{sec:structure} and \ref{sec:map} respectively.

\newpage

\pagebreak

\pagebreak

 \section{Bound on the Weight of the Repair}\label{sec:repair}
 After exposing the minimal outgoing edge from each vertex except $\partial$, we are left with a directed graph that has a very precise structure. Recall, we are denoting this graph by $\Kn^{1}$.
\begin{defn}
A \textbf{cycle-bounded arborescence} is a pair $(\vec M,\vec C)$, where $\vec M$ is a directed graph in which every vertex has exactly one outgoing edge, and $\vec C$ is a directed cycle in $\vec M$ such that the directed graph obtained by contracting all vertices of $\vec C$ to a single vertex $\partial_{\vec C}$ and deleting the resulting self-loops is an arborescence with boundary $\partial_{\vec C}$. 
\end{defn}

\begin{figure}[H]
    \hspace{4.5cm}
    \includegraphics[width=0.9\linewidth]{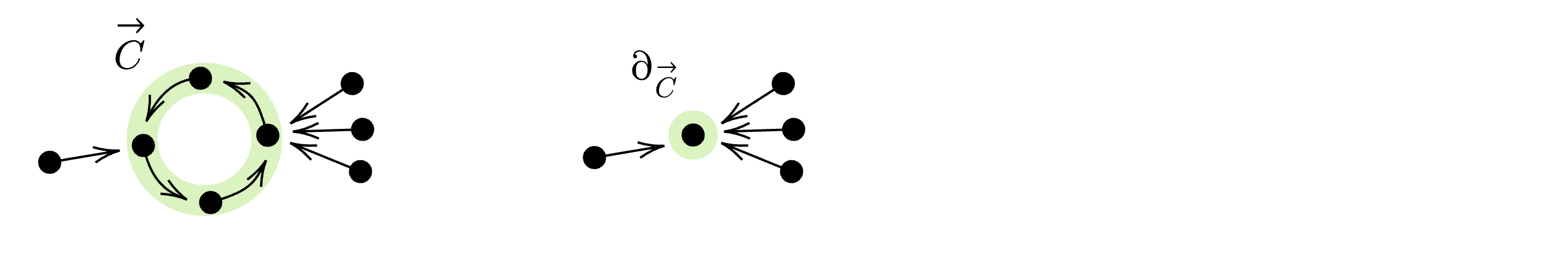}
    \caption{\textbf{Left.} A cycle-bounded arborescence. \textbf{Right.} The arborescence with boundary $\partial_{\vec C}$ obtained by contracting the vertices in its cycle $\vec C$.}
    \label{fig:contraction}
\end{figure}

The connected components of $\vec{K}_{n}^{1}$ can now be completely characterized. See Figure \ref{fig:kn1}. Both of the following facts can be proved using a straightforward induction argument:
\begin{enumerate}
    \item Every component $\cC$ in $\Kn^1$ which is not $\cbd$ is a cycle-bounded arborescence.    
    \item The component $\cbd$ is an arborescence with boundary $\partial$.    
\end{enumerate}

For individual components, the following holds. 
\begin{lemma}\label{lem:compstructure}
Consider the event that a collection of vertices $U=\{v_{1},v_{2},\ldots,v_{m}\}\subseteq V_n$ not containing $\partial$ are the vertex set of  a component $\cC$ of $\Kn^1$. Conditioned on this event, $\cC$ is equal in distribution to a uniform cycle-bounded arborescence on U.  
\end{lemma}
\begin{proof}[Sketch of proof]
Let $\Kn^{\, "}$ denote the graph obtained by independently choosing an outgoing edge uniformly from every vertex except $\partial$. Then $\Kn^{1}$ and $\Kn^{\, "}$ are equal in distribution (note we are looking only at the edges themselves, not the weights). This is because the weights under consideration are i.i.d.\ Exponential$(1)$, and thus continuous,  so any integer linear combination of them will not be zero almost surely (in particular, no two edges will have the same weight).
\end{proof}

More useful structural properties of $\Kn^1$ are proved in Sections \ref{sec:structure} and \ref{sec:map}. In particular, $\Kn^1$ shares a close relationship with the so-called random mapping model, which we detail in Section \ref{sec:map}.\\

Recall that we denote past of a vertex $u$ in $\cmax$ by ${\sf{Past}}_u$. We will also work with fixed $\epsilon,\delta,\eta>0$ satisfying the following properties:
\begin{align}\label{eq:epsdeletaprop}
0<\delta<\epsilon,\text{ } \max \{2(\epsilon+\delta), 6\delta\}<1 \text{ and }\frac{1}{2}<\eta<1-\frac{\epsilon+\delta}{2}.   
\end{align}

\begin{defn}\label{def:A}
For the largest component $\cmax$, and for $\epsilon$ fixed in agreement with \eqref{eq:epsdeletaprop}, define the following random set of vertices 
\[
\cA_{\epsilon}:=\left\{u\in \cmax\setminus \Cmax: \frac{|\cmax|^{\frac{1}{2}+\epsilon}}{|\Cmax|} \le |{\sf{Past}}_u|\le 2\,\frac{|\cmax|^{\frac{1}{2}+\epsilon}}{|\Cmax|}  \right\}.
\]
\end{defn}

Recall that $|\sf{Comp}|$ is the number of components in $\Kn^1$. For $\epsilon$, $\delta$ and $\eta$ as per \eqref{eq:epsdeletaprop}, we will refer to the following event as \textbf{good}:
\begin{align}\label{eq:good}
\cE_{\epsilon,\delta,\eta}:=\left\{|{\sf Comp}|\le\tfrac{11}{\log 2}\log n, \text{ }|\cmax|^{\frac{1}{2}-\delta}< |\Cmax|< |\cmax|^{\frac{1}{2}+\delta},\text{ }|\cA_{\epsilon}|>\left(\frac{\log 2}{11}\frac{n}{\log n}\right)^{\eta}\right\}.
\end{align}

There are two important properties that are guaranteed to hold on the good event that will be used repeatedly:
\begin{enumerate}
\item $|\cmax|>\frac{n\log 2}{11\log n}$, follows immediately using the pigeonhole principle. 
\item $\cmax\neq \cbd$ for $n$ large enough, since we know that on $\cE_{\epsilon, \delta,\eta}$, $|\Cmax|>(\frac{n\log 2}{11\log n})^{\frac{1}{2}-\delta}>1$, final inequality holds so long as $n\log 2>11\log n$. This allows us to ignore the case $\cmax=\cbd$.  
\end{enumerate}
As suggested by the nomenclature, the good event $\cE_{\epsilon,\delta,\eta}$ is typical. 
\begin{lemma}\label{lem:good}
As $n\to \infty$, 
\[\P(\cE_{\epsilon,\delta,\eta}^{c})=O\left( \left(\frac{n}{\log n}\right)^{\max\big\{-\delta,-1+\max\{\epsilon+\delta, 6\delta\}\big\}}\right). \]
\end{lemma}

The proof will be presented in Section \ref{sec:map}. Before we go further, let us state one useful fact about exponential distribution.
\begin{lemma}\label{lem:memoryless}
Let $X_1,\cdots,X_n$ be i.i.d.\ Exponential$(1)$ random variables and $X_J=\min_{1\leq i\leq n}X_{i}$ be the minimum among them. Then conditioned on $J$, the variables $(X_i-X_J)_{i\neq J}$ are i.i.d.\ Exponential$(1)$. 
\end{lemma}
\begin{proof}
Fix \(j\in\{1,\ldots,n\}\). By the memoryless property of the exponential distribution, conditional on \(J=j\) and \(X_J=x\), the random variables \((X_i-X_J)_{i\ne j}\) are independent and each has distribution \(\operatorname{Exponential}(1)\). Their joint distribution does not depend on \(x\). Integrating over \(X_J\), conditional on \(J=j\), shows that the same conditional joint distribution holds given \(J=j\) alone. This proves the result.
\end{proof}

We now accumulate the errors from each of the repair steps. Note, the exposure of the minimal outgoing edge from each vertex results in a random graph $\vec{K}_{n}^{1}$ (we are ignoring the weights themselves, we are only identifying which of the edges are minimal). Knowledge of $\vec{K}_{n}^{1}$ immediately tells us whether the good event $\cE_{\epsilon,\delta,\eta}$ occurs or not. For notational convenience will denote by $\E^{n,\cE_{\epsilon,\delta,\eta}}$ the conditional expectation given $\vec{K}_{n}^{1}$, restricted to $\cE_{\epsilon,\delta,\eta}$: 
\[
\E^{n,\cE_{\epsilon, \delta,\eta}} (X):={\E(X\mathds{1}_{\cE_{\epsilon, \delta, \eta}}\big|\vec{K}_{n}^{1})}=\mathds{1}_{\cE_{\epsilon, \delta, \eta}}\E(X\big|\vec{K}_{n}^{1}). 
\]

\begin{lemma}\label{lem:errbd1}
Let $\vec K_{n}^{2}$ denote the graph obtained after the completion of \textsf{Step 1}. Note, there are only two components: $\cbd$ and $\cmax$. Let $W_{n}^{1}$ denote the weight of the edges so obtained. Then 
\[
0\leq (W_{n}^{1}-W_{n})\mathds1_{\cE_{\epsilon,\delta,\eta}} \le R^{1}_{n}, \text{ where } \E(R^{1}_{n}) =O\left(\frac{(\log n)^2}{n}\right).
\]
\end{lemma}
\begin{proof}
We condition on $\Kn^1$ and assume that $\cE_{\epsilon,\delta,\eta}$ occurs. Let $\cC_i$ be a component of $\Kn^1$ which is neither $\cbd$ nor $\cmax$. Let ${\sf Cycle}_{i}$ be the unique cycle in $\cC_i$, which is of size at least $2$ since our vertices do not have self-loops. Let $E_{\max}(u)$ denote the collection of edges with tail $u\in {\sf Cycle}_i$ and head in $\cmax$, and let $\vec{e}_{*}(u)$ denote the existing outgoing edge from $u$ in $\vec{K}_{n}^{1}$, in particular the minimal edge at $u$. Define \[\vec e_i^{\,\dagger}(u):=\text{argmin}\{Y_{\vec{e}}-Y_{\vec{e}_{*}(u)}:e\in E_{\max}(u)\}\]
and \[
\vec{e}^{\, \dagger}_{i}:=e^{\, \dagger}(u_{i*})\text{ where }u_{i*}:=\text{argmin}\{Y_{\vec{e}^{\, \dagger}(u)}-Y_{\vec{e}_{*}(u)}:u\in {\sf Cycle_{i}}\}.
\]
Define $\delta_{i}:=Y_{\vec e^{\, \dagger}(u_{i*})}-Y_{\vec{e}_{*}(u_{i*})}$. Note, applying Lemma \ref{lem:memoryless}, and using the fact that the minimum of $m$ exponential variables is equal in distribution to an exponential random variable with mean $m^{-1}$, we obtain that $\delta_{i}$ is an exponential random variable with mean given by $|\cmax|^{-1}\cdot |{\sf Cycle}_{i}|^{-1}$.  Since we are on $\cE_{\epsilon,\delta,\eta}$, the pigeonhole principle gives us $|\cmax|>\frac{n\log 2}{11\log n}$.  We obtain \begin{equation}\label{eq:eimax}
|{\sf Cycle}_i|\cdot |\cmax|\geq 2\cdot \frac{\log 2}{11}\frac{n}{\log n}. 
\end{equation}
Thus, $\E^{n,\cE_{\epsilon,\delta,\eta}} (\delta_{i}) \leq \frac{11\log n}{n\cdot 2\log 2 }$. Since there are  at most $\frac{11}{\log 2}\log n$ many components, this procedure has to be repeated at most $\frac{11}{\log 2}\log n $ times, and the error accumulated is the sum of $\delta_{i}$ over all components.  
Thus, 
\[
0\leq (W_{n}^{1}-W_{n})\mathds1_{\cE_{\epsilon,\delta,\eta}}\le R^{1}_{n}:=\left(\sum_{i}\delta_{i}\right)\mathds1_{\cE_{\epsilon,\delta,\eta}}.
\]
Clearly,
\[
\E^{n,\cE_{\epsilon,\delta,\eta}} (R_{n}^{1})=\sum_{i}\E^{n,\cE_{\epsilon,\delta,\eta}}(\delta_{i})\leq \frac{11\log n}{2n\log 2}\frac{11}{\log 2}\log n=O\left(\frac{(\log n)^2}{n}\right). 
\]

Taking expectation again, we conclude. 
\end{proof}
Next, we bound the weight of {\sf Step 2}, in which we grow the component of $\partial$ by attaching a large arborescence to it. 

\begin{lemma}\label{lem:errbd2}
Let $\Kn^{3}$ denote the graph obtained after the completion of \textsf{Step 2}. Let $W_{n}^{2}$ denote the sum of weights of edges in $\Kn^{3}$. Then
\[
0\leq (W_{n}^{2}-W_{n}^{1})\mathds1_{\cE_{\epsilon,\delta,\eta}}\le R_{n}^{2}, \text{ where }\E(R_{n}^{2})=O\left(\left(\frac{n}{\log n}\right)^{-\eta}\right). 
\]
\end{lemma}
\begin{proof}
 We condition on $\vec{K}_{n}^{1}$ and assume that $\cE_{\epsilon,\delta,\eta}$ occurs. For a given vertex $u$, let the existing edge in $\vec{K}_{n}^{1}$ with tail $u$ be denoted by $\vec{e}_{*}(u)$. Consider the collection of all vertices $u\in \cA$, and the collection of all edges from $u$ to $\cbd$ in $\Kn$.  Denote this collection $E_{\partial}(u)$ and define  
\[
\vec{e}^{\,\dagger}(u)=\text{argmin}\{ Y_{\vec e}-Y_{\vec{e}_{*}}:\vec e\in E_{\partial}(u)\}.
\] 
Note, $|E_{\partial}(u)|=|\cbd|$, so by Lemma \ref{lem:memoryless} $\{Y_{\vec{e}^{\,\dagger}(u)}-Y_{\vec{e}_{*}(u)}\}_{u\in \cA_{\epsilon}}$ are a family of i.i.d. exponential random variables with mean $|\cbd|^{-1}$. Next, consider 
\[
u_{*}:=\text{argmin}\{Y_{\vec{e}^{\,\dagger}(u)}-Y_{\vec{e}_{*}(u)}:u\in \cA\}.
\]
Note, as the minimum of $|\cA|$ many independent exponential random variables with mean $|\cbd|^{-1}$, again by Lemma \ref{lem:memoryless}, $Y_{\vec{e}^{\,\dagger}(u_{*})}-Y_{\vec{e}_{*}(u_{*})}$ is exponential with mean $|\cA|^{-1}|\cbd|^{-1}\leq (n\log 2/11\log n)^{-\eta}$. 
Finally, 
 \[
0\leq (W_{n}^{2}-W_{n}^{1})\mathds1_{\cE_{\epsilon,\delta,\eta}}\leq (Y_{\vec{e}^{\,\dagger}(u_{*})}-Y_{\vec{e}_{*}(u_{*})})\mathds1_{\cE_{\epsilon,\delta,\eta}}=:R_{n}^{2}. 
 \]
The final step is identical to the previous lemma. 
\end{proof}

In the following lemma we bound the weight of the last repair step, in which we connect $\Cmax$ to the component of $\partial$, which is now large enough.
\begin{lemma}\label{lem:errbd3}
Let $\Kn^{4}$ denote the arborescence obtained after the completion of \textsf{Step 3}. Then
\[
0<(W_{n}^{3}-W_{n}^{2})\mathds1_{\cE_{\epsilon,\delta,\eta}}\le R_{n}^{3}, \text{ where }\E(R_{n}^{3})=O\left(\left(\frac{n}{\log n}\right)^{-\frac{1}{2}-\epsilon}\right).
\]
\end{lemma}
\begin{proof}
 We condition on $\vec{K}_{n}^{1}$ and assume that $\cE_{\epsilon,\delta,\eta}$ occurs. We take the edge from the vertices of $\Cmax$ to ${\sf{Past}}_{u_{*}}\cup \cbd$ that minimizes the repair cost. By an argument identical to that of Lemma \ref{lem:errbd2}, we find the repair cost is given by an exponential random variable with mean 
\[|{\sf{Past}}_{u_*}\cup\cbd|^{-1}|\Cmax|^{-1}\leq |\cmax|^{-\frac{1}{2}-\epsilon}\leq \left(\frac{\log 2}{11}\cdot \frac{n}{\log n}\right)^{-\frac{1}{2}-\epsilon},\] which follows from Lemma \ref{lem:memoryless}. The remainder of the argument is identical to the previous two lemmas. 
\end{proof}

To extend to the general class of weights as per the criteria defined in Theorem \ref{thm:main}, all we need to do is establish an analogue for Lemma \ref{lem:memoryless}. The steps are then identical. This analogue is presented in Section \ref{sec:conc}, in particular Lemma \ref{lem:genweight}. We are now in a position to prove the CLT, modulo Lemma \ref{lem:good}.
\begin{proof}[Proof of Theorem \ref{thm:main} in the exponential case]
As demonstrated in Section \ref{sec:proofmethod}, proving the CLT boils down to a verification that $\sqrt{n}(W_{n}^{*}-W_{n})$ converges to $0$ in probability. For a fixed $\varepsilon'>0$ let $\cB_n:=\{ \sqrt{n}\left(W_n^* -W_n\right) > \varepsilon'  \}$ be the event whose probability we want to prove goes to 0. Note,
\[
\cB_n = (\cB_n \cap \cE_{\epsilon,\delta,\eta}) \cup (\cB_n \cap \cE_{\epsilon,\delta,\eta}^c),
\]
which implies, by Lemma \ref{lem:good}, 
\begin{align*}
    \P(\cB_n) & = \P(\cB_n \cap \cE_{\epsilon,\delta,\eta}) + \P(\cB_n \cap \cE_{\epsilon,\delta,\eta}^c) \\[0.3cm]
    & \leq \P(\cB_n \cap \cE_{\epsilon,\delta,\eta}) + \P(\cE_{\epsilon,\delta,\eta}^c) \\ &\leq \P(\cB_{n}\cap \cE_{\epsilon,\delta,\eta})+O\left(\left(\tfrac{n}{\log n}\right)^{\max\big\{-\delta,-1+\max\{\epsilon+\delta, 6\delta\}\big\}}\right).
\end{align*}
To show that $\P(\cB_{n}\cap \cE_{\epsilon,\delta,\eta})$ goes to $0$, we showed that $(W_{n}^{*}-W_{n})\mathds1_{\cE_{\epsilon,\delta,\eta}}$ is upper bounded by a random variable $R_{n}$, where $\sqrt{n}R_{n}\to 0$ in probability. Note that, at the end of \textsf{Step 3}, we have an arborescence (see Figure \ref{fig:repair}) whose total weight $W_{n}^{3}$ upper bounds $W_{n}^{*}$. Thus,
\[
0<(W_{n}^{*}-W_{n})\mathds1_{\cE_{\epsilon,\delta,\eta}}\le (W_{n}^{3}-W_{n})\mathds1_{\cE_{\epsilon,\delta,\eta}}\le R_{n}^{1}+R_{n}^{2}+R_{n}^{3}, 
\]
where $R_{n}^{1},R_{n}^{2}$ and $R_{n}^{3}$ are as per Lemmas \ref{lem:errbd1}, \ref{lem:errbd2} and \ref{lem:errbd3} respectively. Define $R_{n}$ to be the sum of the three error terms. By the prior three lemmas,  
\[
\E(\sqrt{n}R_{n}) =O\left(\frac{(\log n)^{\frac{1}{2}+{\max\{\epsilon, \eta-\frac{1}{2}\}}}}{n^{\min\{\epsilon, \eta-\frac{1}{2}\}}}\right)=o(1).
\]
 Notice that the leading order error term comes from the last two steps (Lemmas \ref{lem:errbd2} and \ref{lem:errbd3}). Thus, $\sqrt{n}(W_{n}^{*}-W_{n})\to 0$ in probability. Since $\sqrt{n}(W_{n}-1)\to \cN(0,1)$ in distribution, by Slutsky's theorem $\sqrt{n}(W_{n}^{*}-1)$ converges to the same limit in distribution.
\end{proof}

\section{Structure of the largest component of 
\texorpdfstring{$\Kn^1$}{Directed Complete Graph}}
\label{sec:structure}
 We need to understand some structural properties of $\cmax$. If we condition  on the event that a specific subcollection of vertices $U\subseteq V_n$  constitutes $\cmax$, then $\cmax$ is equal in distribution to a uniformly sampled cycle-bounded arborescence on those vertices. The previous claim follows from Lemma \ref{lem:compstructure}. Going one step further, conditioned on the event that the vertices constituting $\Cmax$ are $\tilde{U}\subseteq U$, then the distribution of the remainder $\cmax \setminus \Cmax$ is uniformly chosen amongst all forests on $U$, rooted at the vertices of $\tilde{U}$. The edges are then oriented towards the roots.
 
 For a given vertex $u$, we now calculate the distribution of the size of its past $|{\sf{Past}}_u|$. We will denote by $\P^{U,\tilde{U}}$ and $\E^{U,\tilde{U}}$ the conditional probability and conditional expectation with respect to the event that the vertex sets of $\cmax$ and $\Cmax$ are $U$ and $\tilde{U}$ respectively. Moreover, let $|U|=m$ and $|\tilde{U}|=k$.
\begin{lemma}\label{lem:tprob1}
 For any $u\in U\setminus \tilde{U}$, 
\[
\P^{U,\tilde{U}}(|{\sf{Past}}_u|=l) = 
\binom{m-k-1}{l-1} \frac{l^{l-2} (m-l)^{m-l-k}}{m^{m-k-1}},
\]
$1\le l\le m-k$.
\end{lemma}
\begin{proof}
The number of rooted forests on $m$ vertices with $k$ designated roots is given by $km^{m-k-1}$, by Cayley's formula.
Next, choose the descendants: ${\sf{Past}}_u$ cannot contain any of the designated roots. Since $u$ is already in ${\sf{Past}}_u$, we must choose $l-1$ additional vertices from the remaining $m-k-1$ available non-roots. There are $\binom{m-k-1}{l-1}$ ways to do this. By Cayley's formula, there are $l^{l-2}$ distinct trees on these $l$ vertices. When we set $u$ as the root of this subtree, all edges are uniquely directed toward $u$. Next, we need to work with the ancestors of $u$. We must have exactly one edge pointing to a vertex outside of ${\sf{Past}}_u$ (there are $m-l$ such vertices). Thus, there are $m-l$ choices for the parent of $u$. Finally, the remaining $m-l$ vertices must form a valid directed forest pointing toward the $k$ roots. No vertex outside ${\sf{Past}}_u$ can be connected to ${\sf{Past}}_u$ (otherwise it would be part of ${\sf{Past}}_u$). Using our forest-counting formula, the number of ways to arrange the remaining vertices is $k(m-l)^{m-l-k-1}$.\qedhere

\end{proof}
A more advanced version of the same lemma will be required for two-point estimates. In particular, the joint distribution of the descendant trees for two distinct vertices in a forest can be explicitly calculated. There are however three sub cases to be considered.  

\begin{lemma}\label{lem:tprob2}
Let $u$ and $v$ be non-root vertices in $U$ (that is, $u,v\in U\setminus\tilde{U}$). Then 
\vspace{0.2cm}
\[
\P^{U,\tilde{U}}\left(|{\sf{Past}}_u|=l_{1},|{\sf{Past}}_v|=l_{2},\cB\right)= \begin{cases} \frac{(m-k-2)!\cdot l_1^{l_1-2} l_2^{l_2-2} (m-l_1-l_2)^{m-l_1-l_2-k+1}}{(l_1-1)!(l_2-1)!(m-k-l_1-l_2)!\cdot m^{m-k-1}} , \cB=\{{\sf{Past}}_u\cap {\sf{Past}}_v=\emptyset\}\\[0.5cm] 
\frac{(m-k-2)!\cdot l_2^{l_2-2} (l_1-l_2)^{l_1-l_2-1} (m-l_1)^{m-l_1-k}}{(l_2-1)!(l_1-l_2-1)!(m-k-l_1)!\cdot m^{m-k-1}}, \cB=\{v\in {\sf{Past}}_u\}\\[0.5cm] 
\frac{(m-k-2)!\cdot l_1^{l_1-2} (l_2-l_1)^{l_2-l_1-1} (m-l_2)^{m-l_2-k}}{(l_1-1)!(l_2-l_1-1)!(m-k-l_2)!\cdot m^{m-k-1}},\cB=\{ u\in {\sf{Past}}_v\}.
\end{cases}
\]
\end{lemma}
\begin{proof}
Note the following trichotomy, we either have ${\sf{Past}}_v\cap {\sf{Past}}_u=\emptyset$,  ${\sf{Past}}_u\subset {\sf{Past}}_v$ if $u\in{\sf{Past}}_v$, or ${\sf{Past}}_v\subset {\sf{Past}}_u$ if $v\in {\sf{Past}}_u$. There are no other possibilities. We also note that we must have $l_{1}+l_{2}\leq m-k$ for the disjoint case, $1\leq l_{2}<l_{1}\leq m-k$ for the case ${\sf{Past}}_v\subset {\sf{Past}}_u$, and the corresponding reversed expression for the final case in our trichotomy.
\begin{itemize}
\item Case 1. We partition the $m-k-2$ available non-root vertices into  ${\sf{Past}}_u$ and ${\sf{Past}}_v$, leaving $m-k-l_1-l_2$ non-roots outside both. For the vertices constituting each descendant tree, there are $\binom{m-k-2}{l_{1}-1}$ and $\binom{m-k-1-l_{1}}{l_{2}-1}$ choices respectively (for vertices other than $u$ and $v$). Multiplying and simplifying yields the prefactor
\[
\frac{(m-k-2)!}{(l_1-1)!(l_2-1)!(m-k-l_1-l_2)!}.
\]
For each of these pasts, there are  $l_1^{l_1-2}$ and $l_2^{l_2-2}$ possibilities respectively. Next, $u$ and $v$ attach themselves to the forest spanned by the remaining $m-l_1-l_2$ vertices, meaning they independently select parents among these, yielding $(m-l_1-l_{2})^{2}$ ways of doing so. Then, apply Cayley's formula to calculate the number of forests spanned by the remainder, which is $k(m-l_{1}-l_{2})^{m-l_{1}-l_{2}-k-1}$. Finally divide by $km^{m-k-1}$, which is the total number of labeled rooted forests on $m$ vertices and $k$ designated roots.  
\item Case 2. We select the $l_2-1$  vertices for ${\sf{Past}}_v$ and the $l_1-l_2-1$ vertices that belong to ${\sf{Past}}_u$ but not ${\sf{Past}}_v$. There are $\binom{m-k-2}{l_{2}-1}$ and $\binom{m-k-l_{2}-1}{l_{1}-l_{2}-1}$ ways of doing so respectively.  The inner subtree ${\sf{Past}}_v$ has $l_2^{l_2-2}$ configurations. The vertex $v$ must then connect to one of the $l_1-l_2$ vertices in the rest of ${\sf{Past}}_u$, with $(l_{1}-l_{2})$ ways of doing so. Finally, $u$ must connect to the outside forest, with $(m-l_{1})$ ways of doing so. Cayley's formula then tells us that the outside forest has $k(m-l_{1})^{m-l_{1}-k-1}$ many configurations to choose from. 
\item Swap the roles of $u$ and $v$, and apply the previous point. \qedhere
\end{itemize}
\end{proof}

We have all the tools required to examine the occurrence of vertices in $\cmax$ with large pasts. Note, conditioned on the event that the vertices of $\cmax$ are a given $U\subset V_n$, with $|U|=m$, and the vertices constituting $\Cmax$ are a given $\tilde{U}$ with $|\tilde{U}|=k$, we can write 
\begin{equation} \label{eq:sizeA}
|\cA|=\sum_{u\in U\setminus \tilde{U}} \mathds{1}\left(\frac{m^{\frac{1}{2}+\epsilon}}{k}\leq|{\sf{Past}}_u|\leq \frac{2m^{\frac{1}{2}+\epsilon}}{k}\right).
\end{equation} 
\begin{lemma}\label{lem:texp}
For $m$ large, and $m^{\frac{1}{2}-\delta}\leq k \leq m^{\frac{1}{2}+\delta}$, $\epsilon$ and $\delta$ as per \eqref{eq:epsdeletaprop},   
\[
\E^{U,\tilde{U}} (|\cA|)=\Theta\left(  \sqrt{k}m^{\frac{3}{4}-\frac{\epsilon}{2}}\right).
\]
\end{lemma}
\begin{proof}
Note that from \eqref{eq:sizeA} we have
\[
\E^{U,\tilde{U}} (|\cA|)=\sum_{u\in U\setminus\tilde{U}}\P^{U,\tilde{U}}\left(\frac{m^{\frac{1}{2}+\epsilon}}{k}\leq|{\sf{Past}}_u|\leq \frac{2m^{\frac{1}{2}+\epsilon}}{k}\right)=(m-k)\P^{U,\tilde{U}}\left(\frac{m^{\frac{1}{2}+\epsilon}}{k}\leq|{\sf{Past}}_u|\leq \frac{2m^{\frac{1}{2}+\epsilon}}{k}\right).
\]
Thus, we need lower and upper bounds of the probability above. Note that 
\[
\P^{U,\tilde{U}}\left(\frac{m^{\frac{1}{2}+\epsilon}}{k}\leq|{\sf{Past}}_u|\leq \frac{2m^{\frac{1}{2}+\epsilon}}{k}\right)=\sum_{\frac{m^{\frac{1}{2}+\epsilon}}{k}\le l\le\frac{2m^{\frac{1}{2}+\epsilon}}{k}}\P^{U,\tilde{U}}(|{\sf{Past}}_u|=l) . 
\]
 If we expand the factorial terms in the formula from Lemma \ref{lem:tprob1}, we can express the probability of the number of vertices in ${\sf{Past}}_u$ being $l$ as 
\[\P^{U,\tilde{U}}(|{\sf{Past}}_u| = l) = \frac{l^{l-2}}{(l-1)!} \prod_{j=1}^{l-1} \left(1 - \frac{k+j}{m}\right) \left(1 - \frac{l}{m}\right)^{m-l-k}.\]
On taking the logarithm, we obtain 
\[
\log\left( \frac{l^{l-2}}{(l-1)!}\right)+\sum_{j=1}^{l-1}\log \left(1-\frac{k+j}{m}\right)+(m-l-k)\log \left(1-\frac{l}{m}\right). 
\]
In the regime under consideration, both $k$ and $l$ are $o(m)$. We can use Taylor's theorem for the second summand to get
\begin{align*}
\sum_{j=1}^{l-1}\log \left(1-\frac{k+j}{m}\right) &=-\sum_{j=1}^{l-1}\frac{k+j}{m}+O\left(\frac{(k+l)^{2}}{m^{2}}\right)\\[0.4cm] &= -\frac{k(l-1)}{m}-\frac{l(l-1)}{2m} +O\left(\frac{l(k+l)^{2}}{m^{2}}\right).
\end{align*}
We can also use Taylor's theorem for the last summand and get 
\[
(m-l-k)\log \left(1-\frac{l}{m}\right)=-\left(l-\frac{l^{2}+kl}{m} \right)+O\left(l^{2}m^{-1}\right).
\]
Adding these two terms, we find that $-\frac{kl}{m}$ gets canceled and some terms can absorbed into the $O$ errors, and we are left with
\[
-l + \frac{k}{m} +O\left(m^{-2}l(k+l)^{2}\right)+O\left(l^{2} m^{-1}\right). 
\]
Notice that since $m^{\frac{1}{2}-\delta} \leq k \leq m^{\frac{1}{2}+\delta}$, then in this regime for $m$ and $l$, the last three of these terms are uniformly $o(1)$ under the hypotheses \eqref{eq:epsdeletaprop}. Applying Stirling's approximation allows us to conclude
\begin{align}\label{eq:tasymp}
\P^{U,\tilde{U}}(|{\sf{Past}}_u| = l)=\frac{l^{l-2} e^{-l}}{(l-1)!}(1+o(1))=\Theta\left(\frac{1}{l^{3/2}}\right),\end{align}
for some implicit constants $C,c>0$. Remember we are considering 
\[
\frac{m^{\frac{1}{2}+\epsilon}}{k}\leq l\leq \frac{2m^{\frac{1}{2}+\epsilon}}{k} 
\]
so for all such $l$, 
\[
\frac{ck^{3/2}}{2^{3/2}m^{3/2(\frac{1}{2}+\epsilon)}}\leq \P^{U,\tilde{U}}(|{\sf{Past}}_u| = l)\leq \frac{Ck^{3/2}}{m^{3/2(\frac{1}{2}+\epsilon)}}
\]
and thus
\[
\bigg(\frac{m^{\frac{1}{2}+\epsilon}}{k}+1\bigg)  \frac{c'k^{3/2}}{m^{3/2(\frac{1}{2}+\epsilon)}}\leq \sum_{ \frac{m^{\frac{1}{2}+\epsilon}}{k}\le l\le\frac{2m^{\frac{1}{2}+\epsilon}}{k}}\P^{U,\tilde{U}}(|{\sf{Past}}_u|=l)\leq \bigg(\frac{m^{\frac{1}{2}+\epsilon}}{k}+1\bigg) \frac{Ck^{3/2}}{m^{3/2(\frac{1}{2}+\epsilon)}}.
\]
Multiplying by $(m-k)$, and remembering that $k\in(m^{\frac{1}{2}-\delta},m^{\frac{1}{2}+\delta})$, we obtain that 
\[
\E^{U,\tilde{U}}(|\cA|)=\Theta\left(\sqrt{k}m^{\frac{3}{4}-\frac{\epsilon}{2}}\right). \qedhere
\]
\end{proof}
In order to get a concentration estimate for $|\cA|$, we also need a variance upper bound. Getting this upper bound requires a decomposition over various cases. To start with, let $\cP=\cA\times \cA$. Note that a particular realization of a random forest induces a (random) partition of the ordered pairs of vertices as 
\[
\cP=\cP_{nest}\cup\cP_{disc}\cup\cP_{self},
\]
where
\[
\cP_{nest}:=\{(u,v):u,v\in \cA,u\in{\sf{Past}}_v\text{ or } v\in{\sf{Past}}_u,u\neq v\},
\]
\[
\cP_{disc}:=\{(u,v):u,v\in \cA, {\sf{Past}}_u\cap{\sf{Past}}_v=\emptyset\},
\]
and 
\[
\cP_{self}:=\{(u,u):u\in \cA\}.
\]
Thus, 
\[
|\cA|^{2}=|\cP_{nest}| +|\cP_{disc}| +|\cA|,
\]
and 
\begin{align}\label{eq:vardecomp}
\E^{U,\tilde{U}} (|\cA|^{2})=\E^{U,\tilde{U}} (|\cP_{nest}|)+\E^{U,\tilde{U}} (|\cP_{disc}|)+\E^{U,\tilde{U}}(|\cA|). 
\end{align}
\begin{lemma}\label{lem:tvar1}
For $m$ large, $m^{\frac{1}{2}-\delta}\le k\le m^{\frac{1}{2}+\delta}$, and $\delta$ as per \eqref{eq:epsdeletaprop},  we have 
\[
\E^{U,\tilde{U}} (|\cP_{nest}|)\leq C''''m ,
\]
for some $C''''>0$.
\end{lemma}
\begin{proof}
We may restrict to the case of ${\sf{Past}}_v\subset {\sf{Past}}_u$, that is $v\in{\sf{Past}}_u$, since twice the expected number of these pairs yields $\E^{U,\tilde U} (|\cP_{nest}|)$. Thus, 
\[
\frac{1}{2}\E^{U,\tilde{U}} (|\cP_{nest}|)= \sum_{u \neq v} \sum_{\frac{2m^{\frac{1}{2}+\epsilon}}{k}\ge l_1> l_2 \ge \frac{m^{\frac{1}{2}+\epsilon}}{k}} \P^{U,\tilde{U}}(v \in {\sf{Past}}_u, |{\sf{Past}}_u|=l_1, |{\sf{Past}}_v|=l_2).
\]
Note that the probability in the sum cannot depend on the specific labels due to symmetry. There are in particular $2\binom{m-k}{2}$ many such pairs. Thus, applying Lemma \ref{lem:tprob2} 
\[
\frac{1}{2}\E^{U,\tilde{U}} (|\cP_{nest}|)= 2\binom{m-k}{2}\cdot\sum_{\frac{2m^{\frac{1}{2}+\epsilon}}{k}\ge l_1> l_2 \ge \frac{m^{\frac{1}{2}+\epsilon}}{k}} \frac{(m-k-2)!\text{ }l_2^{l_2-2} (l_1-l_2)^{l_1-l_2-1} (m-l_1)^{m-l_1-k}}{(l_2-1)!(l_1-l_2-1)!(m-k-l_1)!\text{ }m^{m-k-1}}.
\]
A simpler expression for the sum can be obtained by conditioning further. Note, 
\begin{align*}
 \P^{U,\tilde{U}}\left(v \in {\sf{Past}}_u ,|{\sf{Past}}_u|=l_1, |{\sf{Past}}_v|=l_2\right) = \P^{U,\tilde{U}} \left(|{\sf{Past}}_v|=l_{2},v \in {\sf{Past}}_u\bigg| |{\sf{Past}}_u|=l_1\right)\\ \times \text{ }\P^{U,\tilde{U}}\left(|{\sf{Past}}_u|=l_1\right).
\end{align*}
By Lemma \ref{lem:tprob1},
\[
\P(|{\sf{Past}}_u|=l_1)=\binom{m-k-1}{l_{1}-1} \frac{l_{1}^{l_{1}-2} (m-l_{1})^{m-l_{1}-k}}{m^{m-k-1}}. 
\]
It simplifies things considerably if we rewrite the summand in terms of the conditional probability. In particular, we obtain
\[
\P^{U,\tilde{U}}\left(|{\sf{Past}}_v|=l_2, v \in {\sf{Past}}_u\bigg| |{\sf{Past}}_u|=l_1 \right)=\frac{l_{1}-1}{m-k-1}\binom{l_{1}-2}{l_{2}-1}\frac{l_2^{l_2-2} (l_1-l_2)^{l_1-l_2-1} }{l_{1}^{l_{1}-2}}.
\]
Indeed, after fixing $u$, we can find other $(l_1-1)$ vertices in ${\sf{Past}}_u$ from $|U\setminus\tilde U|-1=(m-k-1)$ vertices. Now, an argument similar to that of Lemma \ref{lem:tprob1} gives the factor $\binom{l_{1}-2}{l_{2}-1}\frac{l_2^{l_2-2} (l_1-l_2)^{l_1-l_2-1} }{l_{1}^{l_{1}-2}}$. By Stirling's approximation, 
\[
\binom{l_1-2}{l_2-1}=\frac{(l_1-2)!}{(l_2-1)!(l_1-l_2-1)!}\le\frac{C}{\sqrt{2\pi}}\frac{(l_1-2)^{l_1-\frac{3}{2}}}{(l_2-1)^{l_2-\frac{1}{2}}(l_1-l_2-1)^{l_1-l_2-\frac{1}{2}}},
\]
where $C>0$. Combining with the terms that come from Cayley's formula, we obtain
\begin{align*}
\P^{U,\tilde{U}} \left(|{\sf{Past}}_v|=l_{2},v \in {\sf{Past}}_u\bigg| |{\sf{Past}}_u|=l_1\right)\le\frac{C}{\sqrt{2\pi}} l_{2}^{-3/2} \left( \frac{l_{1}}{l_{1}-l_{2}} \right)^{1/2}= \frac{C}{\sqrt{2\pi}} l_{2}^{-3/2} \left( 1 - \frac{l_{2}}{l_{1}} \right)^{-1/2}.
\end{align*}
Recall from \eqref{eq:tasymp},
\[\P^{U,\tilde{U}}(|{\sf{Past}}_u| = l_{1})=\frac{l_{1}^{l_{1}-2} e^{-l_{1}}}{(l_{1}-1)!}(1+o(1))\leq \frac{C'}{\sqrt{2\pi }l_1^{3/2}},\]
for some $C'>0$. Combining, we obtain that for some $C''>0$, 
\begin{align*}
\frac{1}{2}\E^{U,\tilde{U}} (|\cP_{nest}|)&\leq (m-k)\sum_{ \frac{m^{\frac{1}{2}+\epsilon}}{k}\le l_1\le \frac{2m^{\frac{1}{2}+\epsilon}}{k}} \sum_{\frac{m^{\frac{1}{2}+\epsilon}}{k}\le l_2<l_1 }\frac{C''}{2\pi} (l_{1}-1)(l_{2} l_{1})^{-3/2}\left( 1 - \frac{l_{2}}{l_{1}} \right)^{-1/2}\\[0.4cm] 
&=(m-k)\sum_{\frac{m^{\frac{1}{2}+\epsilon}}{k}\le l_1\le \frac{2m^{\frac{1}{2}+\epsilon}}{k}}\frac{C''}{2\pi} (l_{1}-1)l_{1}^{-3}l_1\cdot\frac{1}{l_1}\sum_{\frac{m^{\frac{1}{2}+\epsilon}}{k}\le l_2\le l_{1}-1}\left(\frac{l_{2}}{l_{1}}\right)^{-3/2}\left(1-\frac{l_{2}}{l_{1}}\right)^{-1/2}.
\end{align*}
For the interior sum, denoting $l_{2}/l_{1}$ by $\xi$, observe that the sum can be represented as 
\[
\frac{1}{l_{1}}\sum_{\xi\in T} \xi^{-3/2}(1-\xi)^{-1/2},\text{ where }T=\left\{\frac{n}{l_{1}}:n\in \N, \text{ }\frac{m^{\frac{1}{2}+\epsilon}}{k}\leq n\leq l_{1}-1\right\}.
\]
The interior sum is thus a Riemann sum, since $l_1\le\frac{2m^{\frac{1}{2}+\epsilon}}{k}$, approximating the integral 
\[
\int_{\frac{m^{\frac{1}{2}+\epsilon}}{k\cdot l_{1}}}^{1}\xi^{-3/2}(1-\xi)^{-1/2}d\xi\leq \int_{\frac{1}{2}}^{1}\xi^{-3/2}(1-\xi)^{-1/2}d\xi<\infty.  
\]
Thus, the inner sum can be uniformly bounded above by a constant. Overall, we obtain
\[
\frac{1}{2}\E^{U,\tilde{U}} (|\cP_{nest}|)\leq C'''(m-k)\sum_{\frac{m^{\frac{1}{2}+\epsilon}}{k}\le l_1\le\frac{2m^{\frac{1}{2}+\epsilon}}{k}}\frac{1}{l_{1}}\leq C'''' (m-k),
\]
where the last inequality follows from well-known bounds on harmonic sums. 
\end{proof}
Next, we address the issue of the disjoint pasts of the vertices in $\cA$. 
\begin{lemma}\label{lem:tvar2}
For $m$ large, $m^{\frac{1}{2}-\delta}\le k\le m^{\frac{1}{2}+\delta}$, and $\epsilon$ and $\delta$ as per \eqref{eq:epsdeletaprop},
\[
\E^{U,\tilde{U}} (|\cP_{disc}|)=(\E^{U,\tilde{U}} (|\cA|))^{2}+O\left(m^{\max\{1,1+5\delta-\epsilon\}}\right). 
\]
\end{lemma}
\begin{proof}
By Lemma \ref{lem:tprob2} for $L=\frac{m^{\frac{1}{2}+\epsilon}}{k}$,
\begin{align*}
\E^{U,\tilde{U}} (|\cP_{disc}|)&=\sum_{L\leq l_{1},l_{2}\leq 2L}\sum_{u\neq v}\P^{U,\tilde{U}}\left(|{\sf{Past}}_u|=l_{1},|{\sf{Past}}_v|=l_{2},{\sf{Past}}_u\cap {\sf{Past}}_v=\emptyset \right)\\[0.3cm]
&=(m-k)(m-k-1)\sum_{L\leq l_{1},l_{2}\leq 2L}\frac{(m-k-2)!l_1^{l_1-2} l_2^{l_2-2} (m-l_1-l_2)^{m-l_1-l_2-k+1}}{(l_1-1)!(l_2-1)!(m-k-l_1-l_2)!m^{m-k-1}} .
\end{align*}
Next, observe that 
\begin{align*}
(\E^{U,\tilde{U}} (|\cA|))^2&=\sum_{L\leq l_{1},l_{2}\leq 2L}\sum_{u,v}\P^{U,\tilde{U}}(|{\sf{Past}}_u|=l_{1})\,\P^{U,\tilde{U}}(|{\sf{Past}}_v|=l_{2})\\[0.3cm]
&=(m-k)^{2}\sum_{L\leq l_{1},l_{2}\le 2L} \binom{m-k-1}{l_1-1} \binom{m-k-1}{l_2-1}\frac{l_1^{l_1-2}l_2^{l_2-2} (m-l_1)^{m-l_1-k}(m-l_2)^{m-l_2-k}}{m^{2m-2k-2}}.  
\end{align*}
Define
\[
R(l_{1},l_{2})=\frac{(m-k)(m-k-1)\frac{(m-k-2)!l_1^{l_1-2} l_2^{l_2-2} (m-l_1-l_2)^{m-l_1-l_2-k+1}}{(l_1-1)!(l_2-1)!(m-k-l_1-l_2)!m^{m-k-1}}}{(m-k)^2 \binom{m-k-1}{l_1-1} \binom{m-k-1}{l_2-1}\frac{l_1^{l_1-2}l_2^{l_2-2} (m-l_1)^{m-l_1-k}(m-l_2)^{m-l_2-k}}{m^{2m-2k-2}}},
\]
the ratio of the respective summands. Denote by 
\[
\Xi(l_{1},l_{2}):=\binom{m-k-1}{l_1-1} \binom{m-k-1}{l_2-1}\frac{l_1^{l_1-2}l_2^{l_2-2} (m-l_1)^{m-l_1-k}(m-l_2)^{m-l_2-k}}{m^{2m-2k-2}},
\]
so that 
\[
(\E^{U,\tilde{U}}|\cA|)^{2}=(m-k)^{2}\sum_{L\leq l_{1},l_{2}\leq 2L}\Xi(l_{1},l_{2}). 
\]
The comparison between $(\E^{U,\tilde{U}}|\cA|)^{2}$ and $\E^{U,\tilde{U}} (|\cP_{disc}|)$ can now be set up. Note, 
\begin{align}\label{eq:ratiodef}
\E^{U,\tilde{U}} (|\cP_{disc}|)= (m-k)^{2}\sum_{L\leq l_{1},l_{2}\leq 2L}\Xi(l_{1},l_{2})\cdot R(l_{1},l_{2}). 
\end{align}
We aim to show that $R$ is close to $1$. Simplifying the expression for $R$ by cancelling terms that appear in common, we find that  
\[
R(l_{1},l_{2})= \frac{m^{m-k-1}}{(m-k)!}\frac{(m-l_{1}-l_{2})^{m-l_{1}-l_{2}-k+1}(m-k-l_{1})!(m-k-l_{2})!}{(m-l_{1})^{m-l_{1}-k}(m-l_{2})^{m-l_{2}-k}(m-k-l_{1}-l_{2})!}.
\]
The ratio term can further be split into a product of two terms, composed of the factorial parts and the power parts respectively. That is, define 
\[
F(l_{1},l_{2}):=\frac{(m-k-l_{1})!(m-k-l_{2})!}{(m-k)!(m-k-l_{1}-l_{2})!}
\] and 
\[
P(l_{1},l_{2})= \frac{m^{m-k-1}(m-l_{1}-l_{2})^{m-l_{1}-l_{2}-k+1}}{(m-l_{1})^{m-l_{1}-k}(m-l_{2})^{m-l_{2}-k}}. 
\]
For ease of notation, denote $m-k$ by $q$. By Stirling's approximation, 
\begin{align*}
F(l_1,l_2)&=\frac{(q-l_{1})^{q-l_{1}}(q-l_{2})^{q-l_{2}}e^{-2q+l_{1}+l_{2}}\sqrt{2\pi(q-l_{1})(q-l_{2})}}{q^q (q-l_{1}-l_{2})^{q-l_{1}-l_{2}}e^{-2q+l_{1}+l_{2}}\sqrt{2\pi q(q-l_{1}-l_{2})}}(1+O(m^{-1}))\\[0.3cm]
&=\frac{(q-l_{1})^{q-l_{1}}(q-l_{2})^{q-l_{2}}}{q^q (q-l_{1}-l_{2})^{q-l_{1}-l_{2}}}(1+O(L\cdot m^{-1}))  
\end{align*}
For $P$, we have 
\begin{align*}
P(l_{1},l_{2})&=\frac{(q+k)^{q-1}(q+k-l_{1}-l_{2})^{q-l_{1}-l_{2}+1}}{(q+k-l_{1})^{q-l_{1}}(q+k-l_{2})^{q-l_{2}}}\\[0.3cm] &=\frac{q+k-l_{1}-l_{2}}{q+k}\cdot \frac{(q+k)^{q}(q+k-l_{1}-l_{2})^{q-l_{1}-l_{2}}}{(q+k-l_{1})^{q-l_{1}}(q+k-l_{2})^{q-l_{2}}} \\[0.3cm]
&=(1+O(L\cdot m^{-1}))\cdot \frac{(q+k)^{q}(q+k-l_{1}-l_{2})^{q-l_{1}-l_{2}}}{(q+k-l_{1})^{q-l_{1}}(q+k-l_{2})^{q-l_{2}}}. 
\end{align*}
Combining terms, we obtain 
\[
R(l_{1},l_{2})=\left( \frac{q-l_{1}}{q+k-l_{1}}\right)^{q-l_{1}}
\left( \frac{q-l_{2}}{q+k-l_{2}}\right)^{q-l_{2}}
\left(\frac{q+k}{q}\right)^{q}
\left( \frac{q+k-l_{1}-l_{2}}{q-l_{1}-l_{2}}\right)^{q-l_{1}-l_{2}}
(1+O(L\cdot m^{-1})).
\]
Observe that 
\begin{align*}
\left( \frac{q-l_{1}}{q+k-l_{1}}\right)^{q-l_{1}}&=\exp \left(-(q-l_{1})\log \left(1+\frac{k}{q-l_{1}}\right)\right)\\[0.3cm] &=\exp\left(-k+\frac{k^{2}}{2(q-l_{1})}-\frac{k^{3}}{3(q-l_{1})^{2}}+O\left(\frac{k^{4}}{q^{3}}\right)\right)
\end{align*}
Applying the same principle for all four terms, we gather our errors. Note, the $k$ terms cancel. Next, the $k^{2}$ terms in the exponential can be combined as 
\[
\frac{k^{2}}{2}\left(\frac{1}{q-l_{1}}+\frac{1}{q-l_{2}}-\frac{1}{q}-\frac{1}{q-l_{1}-l_{2}}\right)=\frac{-k^{2}l_1l_2(2q - l_1 - l_2)}{2q(q-l_1)(q-l_2)(q-l_1-l_2)}.
\]
Finally we deal with the $k^{3}$ terms for which the relevant error contribution is 
\[
\frac{k^{3}}{3}\left(-\frac{1}{(q-l_{1})^{2}}-\frac{1}{(q-l_{2})^{2}}+\frac{1}{q^{2}}+\frac{1}{(q-l_{1}-l_{2})^{2}}\right).
\]
The expanded term is more complicated, and is given by 
\[
\frac{k^{3}l_1l_2 \left( 6q^4 - 12(l_1+l_2)q^3 + (8l_1^2 + 18l_1l_2 + 8l_2^2)q^2 - 2(l_1+l_2)(l_1^2 + 3l_1l_2 + l_2^2)q + l_1l_2(l_1+l_2)^2 \right)}{3q^2(q-l_1)^2(q-l_2)^2(q-l_1-l_2)^2}.
\]
 Plugging in the fact that $l_{1}$ and $l_{2}$ are $\Theta(m^{\frac{1}{2}+\epsilon}/k)$, and recalling that the regime of $k$ in question is $m^{1/2-\delta}\le k\le m^{1/2+\delta}$ and keeping the dominant errors, we obtain  
\begin{align*}
R(l_{1},l_{2})=\exp \left( O\left(\frac{1}{m^{2-2\epsilon}}\right)+O\left(\frac{k}{m^{3-2\epsilon}}\right)+O\left(\frac{k^{4}}{m^{3}}\right)\right)\left(1+O\left(\frac{1}{km^{\frac{1}{2}-\epsilon}}\right)\right). 
\end{align*}
Note that the above expression is uniform in $(l_{1},l_{2})$. For the range of $k$ in question, the $k^{4}/m^{3}$ term dominates. We may now recast \eqref{eq:ratiodef} as 
\begin{align*}
\E^{U,\tilde{U}} (|\cP_{disc}|) & = (m-k)^{2}\sum_{L\leq l_{1},l_{2}\leq 2L}\Xi(l_{1},l_{2})\cdot R(l_{1},l_{2})\\[0.3cm] 
&\left(1+O\left(\frac{k^{4}}{m^{3}}\right)+O\left(\frac{1}{km^{\frac{1}{2}-\epsilon}}\right)\right)\cdot (m-k)^{2}\sum_{L\leq l_{1},l_{2}\leq 2L}\Xi(l_{1},l_{2})\\[0.3cm]
& =\left(1+O\left(\frac{k^{4}}{m^{3}}\right)+O\left(\frac{1}{km^{\frac{1}{2}-\epsilon}}\right)\right)(\E^{U,\tilde{U}}(|\cA|))^{2} \\[0.3cm]
& = \left(1+O\left(\frac{k^{4}}{m^{3}}\right)+O\left(\frac{1}{km^{\frac{1}{2}-\epsilon}}\right)\right)(\E^{U,\tilde{U}}(|\cA|))^{2}
\end{align*}
 Since $(\E^{U,\tilde{U}}(|\cA|))^{2}=O(m^{\frac{3}{2}-\epsilon}k)$ by Lemma \ref{lem:texp}, 
we obtain 
\[
\E^{U,\tilde{U}} (|\cP_{disc}|)=(\E^{U,\tilde{U}} (|\cA|))^{2}+O\left(\frac{k^{5}}{m^{\frac{3}{2}+\epsilon}}\right)+O(m)=(\E^{U,\tilde{U}} (|\cA|))^{2}+O\left(m^{\max\{1,1+5\delta-\epsilon\}}\right). 
\qedhere
\]
 
\end{proof}
\begin{corollary}
For $m$ large, $m^{\frac{1}{2}-\delta}\le k\le m^{\frac{1}{2}+\delta}$, and $\epsilon$ and $\delta$ as per \eqref{eq:epsdeletaprop}, 
\[
\var^{U,\tilde{U}} (|\cA|)\leq Cm^{\max\{1,1+5\delta-\epsilon\}},
\]
for some $C>0$.
\end{corollary}
\begin{proof}
Note, 
\[
(\E^{U,\tilde{U}}(|\cA|^{{2}}))=\E^{U,\tilde{U}} (|\cP_{nest}|)+\E^{U,\tilde{U}} (|\cP_{disc}|)+\E^{U,\tilde{U}}(|\cA|). 
\]
by \eqref{eq:vardecomp}. Applying Lemmas \ref{lem:tvar1}, \ref{lem:tvar2} and \ref{lem:texp}, we obtain 
\[
\E^{U,\tilde{U}}(|\cA|^{2})\leq Am +(\E^{U,\tilde{U}}(|\cA|))^{2}+A'm^{\max \{1,1+5\delta-\epsilon\}}+A''m^{3/4-\frac{\epsilon}{2}}\sqrt{k}.
\]
Subtracting $(\E^{U,\tilde{U}} (|\cA|))^{2}$ from both sides completes the proof.  
\end{proof}
\begin{corollary}\label{cor:Abound}
As $m$ becomes large, there exists $C'<\infty$ such that 
\[
\P^{U,\tilde{U}}(|\cA|\leq m^{\eta})\leq  C' m^{\max\{-1+\epsilon+\delta,-1+6\delta\}}.
\]

\end{corollary}
    \begin{proof}
    By Chebyshev's inequality,
    \[
    \P^{U,\tilde{U}}(|\cA| \leq m^{\eta}) \le \P^{U,\tilde{U}}(\big||\cA| - \E^{U,\tilde{U}}(|\cA|)\big| \geq \E^{U,\tilde{U}}(|\cA|) - m^{\eta}) \leq c \frac{\text{Var}^{U,\tilde{U}}(|\cA|)}{\left(\E^{U,\tilde{U}}(|\cA|)\right)^2},\]
    for some $c>0$. This is because by our choice of parameters \eqref{eq:epsdeletaprop}, $\eta<1-\frac{\epsilon+\delta}{2}$, and thus $\E^{U,\tilde{U}} (|\cA|)$ dominates $m^{\eta}$ in the limit. Plugging in our bounds, we obtain
    \[
    \frac{\text{Var}^{U,\tilde{U}}(|\cA|)}{\left(\E^{U,\tilde{U}}(|\cA|)\right)^{2}} \leq C' \frac{  m^{\max \{1,1+5\delta-\epsilon\}}}{  m^{3/2 - \epsilon} k } = C'm^{\max\{-\frac{1}{2}+ \epsilon,-\frac{1}{2}+5\delta\}} k^{-1}\leq C' m^{\max\{-1+\epsilon+\delta,-1+6\delta\}}. \qedhere
    \]
    \end{proof}

\section{The Random Mapping Model and Typicality of
\texorpdfstring{$\cE_{\epsilon,\delta,\eta}$}{Good Event}}
\label{sec:map}
Let $(K_{n}, \partial)$ denote the undirected complete graph on $n$ vertices with a distinguished vertex (root) $\partial$. We define the random mapping to be a random directed graph obtained as follows: define $\cM$ to be the collection of all functions $f:V_{n} \to V_{n}$ such that $f(\partial)=\partial$. Note, there are $n^{n-1}$ many such distinct functions. For  $f \in \cM$, we associate a directed graph, by placing the directed edge $(u,v)$ vertices $u$ and $v$ if and only if $f(u)=v$. We ignore the self-loop at $\partial$. Since we are working with functions, note that every non-root vertex has an out-degree of exactly $1$. We will denote this directed graph by $\Kn^{map}$.\\

 The following are known results on random mappings.
\begin{lemma}
We say that a vertex $v$ is cyclic for $F$ if $F^{(k)}(v)=v$ for some $k\in \N$. Let $L_{n}$ denote the number of cyclic vertices associated to a uniformly sampled $F$. We have that 
\[
\P(L_{n}=k)=\frac{kn^{n-k-1}\binom{n-1}{k-1}(k-1)!}{n^{n-1}}=\frac{k}{n}\cdot \prod_{j=1}^{k-1}\left(1-\frac{j}{n}\right).
\]
\end{lemma}
\begin{proof}
There are $\binom{n-1}{k-1}$ different choices for the $k-1$ cyclic vertices, since $\partial$ is already counted. Cyclic vertices must get mapped to each other, and moreover the mappings must be one to one. Thus $F$ acts as a permutation on the cyclic vertices other than $\partial$, and there are $(k-1)!$ choices for the permutation. We now need to consider the possibilities for the other vertices. This corresponds to a labeled spanning forest, where the cyclic vertices are roots. There are $kn^{n-k-1}$ such forests.
\end{proof}
This gives us a robust bound on the typical number of cyclic vertices. 
\begin{lemma}\label{lem:numberofcyclic}
For $a\leq \frac{1}{2}-\delta $, \[\P(L_{n}< n^{a})\leq n^{-2\delta}.  \]
\end{lemma}
\begin{proof}
We know that 
\begin{align*}
\P(L_{n}< n^{a})&=\sum_{1\le k< n^{a}}\frac{k}{n}\cdot \prod_{j=1}^{k-1}\left(1-\frac{j}{n}\right)\\[0.3cm] &\leq \sum_{1\le k< n^{a}}\frac{k}{n} \\[0.3cm] &\leq \frac{1}{n}(n^a-1)(n^{a}-2)\leq n^{2a-1}\leq n^{-2\delta}. \qedhere  
\end{align*}

\end{proof}
Bounding the number of cyclic vertices allows us to control the number of connected components of the random mapping model, and thus our model as well. 
\begin{lemma}\label{lem:compnumber}
Let ${\sf Comp}$ denote the collection of connected components of $\vec{K}_{n}^{map}$. Then as $n\to \infty$,  
\[
\P\left(\frac{{\sf |Comp|}}{\log n}>\frac{11}{\log 2}\right)\leq n^{-2\delta}.
\]
\end{lemma}
This result is rather well known, so we omit a detailed proof and provide a sketch. Later, we will argue that the same bound applies to $\sf Comp$ as defined for $\vec{K}_{n}^{1}$. Thanks to a coupling commonly associated to Feller, appearing in \cite{Feller1945}, it is known that the number of cycles of a uniform permutation on $m$ letters can be represented as a sum of independent Bernoulli random variables. In particular, if $Z_{m}$ denotes the number of cycles, then 
\[
Z_{m}=\sum_{k=1}^{m} I_{k}, \text{ }I_{k}\sim\text{Ber}(k^{-1}). 
\]
In particular, this allows us to conclude two things, the first being that 
\[
\E(Z_{m})=\log m +O(1), 
\]
by the well known asymptotics of the harmonic sum. Moreover, observe that as a consequence of the independence,  
\[
\E \left(2^{Z_{m}}\right)=\prod_{k=1}^{m}\left(\frac{k-1}{k}+\frac{2}{k}\right)=\prod_{k=1}^{m}\left(\frac{k+1}{k}\right)=m+1.
\]
By Markov's inequality,
\[
\P(Z_{m}\geq t \log m)\leq \frac{m+1}{2^{t\log m}}. 
\]
Thus, taking $t=11/\log 2$, we find that 
\[
\P\left(\frac{Z}{\log m}\geq \frac{11}{\log 2}\right)\leq m^{-10}.
\]
Recall, the restriction of a random mapping to the cyclic vertices is a uniform permutation on those vertices. Since the cyclic vertices act as roots for the other vertices, the disjoint cycles of the permutations also index the connected components of the graph corresponding to the random map.  Thus, given that the number of cyclic vertices is $L_{n}=m$, we immediately have that 
\[
\P\left(|{\sf Comp}|> \frac{11}{\log 2}\log m\bigg|L_{n}=m\right)\leq m^{-10}. 
\]
Finally, observe that 
\begin{align*}
\P\left(|{\sf Comp}|\leq \frac{11}{\log 2}\log n\right)&\geq \sum_{m\ge n^{\frac{1}{2}-\delta}}\P\left(|{\sf Comp}|\leq  \frac{11}{\log 2}\log m\bigg|L_{n}=m\right)\;\P(L_{n}=m)\\[0.4cm] &\geq 
\sum_{m\geq n^{\frac{1}{2}-\delta}}(1-n^{-5+10\delta})\;\P(L_{n}=m)\\[0.4cm] &\geq (1-n^{-5+10\delta})(1-n^{-2\delta})
\end{align*}
The final step follows from Lemma \ref{lem:numberofcyclic}. An important corollary is the following.
\begin{corollary} \label{cor:maxsize}
As $n\to \infty$,
\[
\P\left(|\cmax|>\frac{n\log 2}{11\log n}\right)\ge 1-O\left(n^{-2\delta}\right).
\]

\end{corollary}
\begin{proof}
This follows from a direct application of Lemma \ref{lem:compnumber} and the pigeonhole principle.  
\end{proof}

Given a component of size $m$, we now need to understand the size of the embedded cycle in that component. The procedure is almost identical to understanding the number of cyclic vertices. We now condition on the event that a collection of vertices $U$ form the vertex set of a non boundary component. On this event, the structure of $\vec{K}_{n}^{1}$ restricted to these vertices is uniform amongst all cycle-bounded arborescences on those vertices, where the cycle length $k$ must be greater than or equal to $2$. This is the only difference as compared to the $\vec{K}^{1}_{n}$ case, the cycle can be a self loop corresponding to $k=1$. We denote this conditional measure as $\P^{U}$, just as in Section \ref{sec:structure}.
\begin{lemma}\label{lem:Csize}
Let $C_{m}$ denote the size of the cycle in a uniformly sampled cycle-bounded arborescence $U$ where $|U|=m$, and self loops are not allowed. Then 
\[
\P^{U}(C_{m}=k)=\frac{\binom{m}{k} k!\cdot m^{m-k-1}}{\sum_{i=2}^{m}\binom{m}{i} i!\cdot m^{m-i-1}}.
\]
\end{lemma}
\begin{proof}
First, we count the possibilities corresponding to $C_{m}=k$. Note that there are $\binom{m}{k}$ choices for the $k$ cyclic vertices. Next, there are $(k-1)!$ many ways of ordering the cycle, up to shifts. Finally, by Cayley's formula, there are $km^{m-k-1}$ ways of constructing labeled forest with our cyclic vertices as roots. This yields
\[
\binom{m}{k} (k-1)!\cdot k\cdot m^{m-k-1}.
\]
Next, the total number of cycle-bounded arborescences can now be obtained by summing the above expression over all $k$.  
\end{proof}
With this in mind, we can now compute a bound on the probability of a cycle being small. 
\begin{lemma}\label{lem:cyctail}
As $m$ becomes large, 
\[
\P^{U}\big(C_{m}\notin (m^{\frac{1}{2}-\delta},m^{\frac{1}{2}+\delta})\big)\leq  m^{-\delta}.
\]
\end{lemma}
\begin{proof}
Note that 
\[
\P^{U}(C_{m}\leq m^{\frac{1}{2}-\delta})=\sum_{2\le k\le m^{\frac{1}{2}-\delta}}\frac{\binom{m}{k} k!\cdot m^{m-k-1}}{\sum_{i=2}^{m}\binom{m}{i} i!\cdot m^{m-i-1}}=\sum_{2\le k\le m^{\frac{1}{2}-\delta}}\frac{\prod_{j=2}^{k-1}\left(1-\frac{j}{m}\right)}{\sum_{i=2}^{m}\binom{m}{i} i!\cdot m^{-i}}.
\]
The trivial upper bound for the numerator suffices, 
\begin{align}\label{eq:cycnum1}
\sum_{2\le k\le m^{\frac{1}{2}-\delta}}\prod_{j=2}^{k-1}\left(1-\frac{j}{m}\right)\leq m^{\frac{1}{2}-\delta}.
\end{align}
For $C_{m}\geq m^{\frac{1}{2}+\delta}$, we need the asymptotics of the numerator where the sum ranges over $k\in (m^{\frac{1}{2}+\delta},m]$. In this case, note that 
\[
\sum_{m^{\frac{1}{2}+\delta}\le k\le m}\prod_{j=2}^{k-1}\left(1-\frac{j}{m}\right)\leq \sum_{m^{\frac{1}{2}+\delta}\le k\le m}\prod_{2\le j\le \min\{k-1,m^{3/5}\}}\left(1-\frac{j}{m}\right).
\]
The choice of $3/5$ as the exponent of truncation is made for a purely aesthetic reason, it is the nicest fraction between $1/2$ and $2/3$. Next, observe that we have 
\begin{multline*}
\sum_{m^{\frac{1}{2}+\delta}\le k\le m}\prod_{2\le j\le\min\{k-1,m^{3/5}\}}\left(1-\frac{j}{m}\right)\\\leq \sum_{m^{\frac{1}{2}+\delta}\le k\le m}\exp\left(-\sum_{2\le j\le\min\{k-1,m^{3/5}\}}\frac{j}{m} +O\bigg(\frac{({\min\{k-1,m^{3/5}\})^{3}}}{m^{2}}\bigg)\right).
\end{multline*}
The smallest value of $k$ in the outer sum is $m^{\frac{1}{2}+\delta}$, and this yields an upper bound (since the sign in the exponential ahead of the dominant term is negative, the error vanishes since $m^{9/5}<m^{2}$). In particular, we obtain 
\begin{align}\label{eq:cycnum2}
\sum_{m^{\frac{1}{2}+\delta}\le k\le m}\prod_{j=2}^{k-1}\left(1-\frac{j}{m}\right)\leq m\cdot \exp\left(-\frac{1}{2}m^{\min\{2\delta,\frac{1}{5}\}}\right).
\end{align}
The last inequality holds for $m$ sufficiently large. We then need a lower bound on the denominator, which is given by 
\[
\sum_{i=2}^{m}\frac{m!}{(m-i)!}m^{-i}=\sum_{i=1}^{m}\frac{m!}{(m-i)!}m^{-i}-1=:I(m)-1.
\]
The sum can be represented by an integral, which will enable the use of Laplace asymptotics. Observe that  
\[
\int_{0}^{\infty} e^{-x} \left(1+\frac{x}{m}\right)^{m}dx=\int_{0}^{\infty}e^{-x}\sum_{i=0}^{m}\binom{m}{i} \left(\frac{x}{m}\right)^{i}dx=I(m)+1. 
\]
Now, note that 
\[
\int_{0}^{\infty} e^{-x} \left(1+\frac{x}{m}\right)^{m}dx=m\int_{0}^{\infty}e^{-mu}(1+u)^{m}du=m\int_{0}^{\infty}e^{-mg(u)}du,
\]
where $g(u)=u-\log (1+u)$. Note that $g$ is increasing, $g(0)=0$, $g'(0)=0$ and is convex. Moreover, $g''(0)=1>0$. Thus, 
\[
\int_{0}^{\infty} e^{-mg(u)}du\geq \frac{C}{\sqrt{m}} 
\]
where the constant $C$ depends on $g''(0)>0$. Thus, 
\begin{align}\label{eq:cycden}
I(m)\geq C\sqrt{m}-1.
\end{align} 
For the lower tail we  obtain by combining \eqref{eq:cycnum1} and \eqref{eq:cycden}
\[
\P^{U}(C_{m}\leq m^{\frac{1}{2}-\delta})\leq \frac{m^{\frac{1}{2}-\delta}}{C\sqrt{m}-2}\leq \frac{m^{-\delta}}{C'}. 
\]
For the upper tail, we obtain
\[
\P^{U}(C_{m}\geq m^{\frac{1}{2}+\delta})\leq m\cdot\exp\left(-\frac{1}{2}m^{\min\{2\delta,\frac{1}{5}\}}\right)\leq m^{-\delta},
\]
the final inequality holding for $m$ sufficiently large. 
\end{proof}
We now have all the ingredients required to show that our good event $\cE_{\epsilon,\delta,\eta}$ is typical for $\vec{K}_{n}^{1}$. The main difference between $\Kn^{map}$ and $\Kn^{1}$ is that self-loops are permitted for the random mapping, corresponding to the fixed points of the associated functions. To set up a correspondence, let 
\[
\cM_{0}=\{f:V_n\to V_n: f(\partial)=\partial, f(u)\neq u\text{ }\forall u\in V_{n}\setminus{\{\partial\}}\}.
\]
Note that for every realization of $\Kn^1$, we can associate a unique random mapping in $\cM_{0}$. Thus, the graph associated to the uniform measure on $\cM_{0}$ is equal in distribution to $\Kn^1$, with the self-loop at $\partial$ ignored. It is also relatively straightforward to see that 
\[
|\cM_{0}|=(n-1)^{n-1}.
\]
The following lemma is immediate:
\begin{lemma}
Let $\cM$ denote the set of all functions $f:V_{n}\to V_{n}$ such that $f(\partial)=\partial$. Let $F$ be uniformly sampled from $\cM$. Then 
\[
\P(F\in \cM_{0}) = \frac{(n-1)^{n-1}}{n^{n-1}} \searrow e^{-1}, 
\]
as $n\to \infty$. 

\end{lemma}
This lemma enables us to translate results from the random mapping model to our setting, especially probability bounds for rare events. The following is immediately clear for a random mapping $F$ and a collection $\cA$, 
\begin{align}\label{eq:maptotree}
\P(F\in\cA|F\in \cM_{0})\leq e\cdot \P(F\in \cA).
\end{align}
In particular, Lemma \ref{lem:compnumber}  and any corollaries all immediately admit corresponding versions for $\vec{K}_{n}^{1}$. When we refer to these lemmas in the following proof, we are referring to the versions for $\vec{K}_{n}^{1}$, in combination with \eqref{eq:maptotree}.
\begin{proof}[Proof of Lemma \ref{lem:good}]
Denote by 
\begin{align*}
\cE_{1}=\left\{|{\sf Comp}|\le\frac{11}{\log 2}\log n\right\}&,\text{ }\cE_{2}=\{|\cmax|^{\frac{1}{2}-\delta}< |\Cmax|< |\cmax|^{\frac{1}{2}+\delta}\}\\[0.2cm] &\text{ and }\cE_{3}=\left\{|\cA_{\epsilon}|>\left(\frac{\log 2}{11}\frac{n}{\log n}\right)^{\eta}\right\}.
\end{align*}
Clearly, 
\[
\cE_{\epsilon,\delta,\eta}=\cE_{1}\cap \cE_{2}\cap \cE_{3},
\]
and thus 
\[
\cE_{\epsilon,\delta,\eta}^{c}=\cE_{1}^{c}\cup(\cE_{2}^{c}\cap \cE_{1})\cup (\cE_{3}^{c}\cap \cE_{2}\cap\cE_{1}).
\]
Applying Lemma \ref{lem:compnumber}, we know that 
\[
\P(\cE_{1}^{c})=O(n^{-2\delta}).
\]
Next, for $\P(\cE_{2}^{c}\cap \cE_{1})$, we split the probability of the intersection in terms of conditioning, 
\begin{align*}
\P(\cE_{2}^{c}\cap \cE_{1})
\leq\sum_{m\geq \frac{n}{t\log n}}\sum_{|U|=m} \P(\cE_{2}^{c}|\cmax=U)\P(\cmax=U),
\end{align*}
where $t=11/\log 2$. We begin the sum at $\frac{n}{t\log n}$ because of the pigeonhole argument of Corollary \ref{cor:maxsize}. By Lemma \ref{lem:cyctail}, on the event $|U|=m$,
\[
\P(\Cmax\notin  (m^{\frac{1}{2}-\delta},m^{\frac{1}{2}+\delta})|\cmax=U)\leq Cm^{-\delta} .
\]
Thus, 
\begin{align*}
\P(\cE_{2}^{c}\cap \cE_{1})&\leq Ct^{-\delta}\left(\frac{n}{\log n}\right)^{-\delta}\sum_{m\geq \frac{n}{t\log n}} \sum_{|U|=m}\P(\cmax=U)\\[0.5cm]
&=O\left(\left(\frac{n}{\log n}\right)^{-\delta}\right). 
\end{align*}

Finally, we need to deal with $\P(\cE_{3}^{c}\cap \cE_{2}\cap \cE_{1})$. For ease of notation, we define 
\[
\cE_{m,k}:=\{|\cmax|=m,|\Cmax|=k\}.
\]
Note that we may further disjointly partition $\cE_{m,k}$ as 
\[
\cE_{m,k}=\bigcup_{\substack{|U|=m \\ \tilde{U}\subset U,|\tilde{U}|=k}}\{\cmax=U, \Cmax=\tilde{U}\}. 
\]

Note that 
\begin{align*}
\P\left(\cE_{3}^{c}\cap \cE_{1}\cap \cE_{2}\right) & \leq \sum_{m\geq \frac{n}{t\log n}}\sum_{m^{\frac{1}{2}-\delta}\le k\le m^{\frac{1}{2}+\delta}}\sum_{\substack{|U|=m \\ \tilde{U}\subset U,|\tilde{U}|=k}}\P(\{\cmax=U\}\cap \{\Cmax=\tilde{U}\}\cap \{|\cA_{\epsilon}|\leq m^{\eta}\})\\[0.3cm]
& =\sum_{m\geq \frac{n}{t\log n}}\sum_{m^{\frac{1}{2}-\delta}\le k\le m^{\frac{1}{2}+\delta}}\sum_{\substack{|U|=m \\ \tilde{U}\subset U,|\tilde{U}|=k}}\P^{U,\tilde{U}}( |\cA_{\epsilon}|\leq m^{\eta}) \; \P(\cmax=U,\Cmax=\tilde{U}).
\end{align*}
Next, we apply Corollary \ref{cor:Abound}. Plugging in, we obtain
\begin{align*}
\P\left(\left\{|\cA_{\epsilon}|\le \left(\frac{n}{t\log n}\right)^{\eta}\right\}\cap\cE_1\cap\cE_2\right) &\leq \sum_{m\geq \frac{n}{t\log n}}\sum_{k=\lceil m^{\frac{1}{2}-\delta}\rceil}^{\lfloor m^{\frac{1}{2}+\delta}\rfloor} Cm^{-1+{\max\{\epsilon+\delta,6\delta\}}} \sum_{\substack{|U|=m \\
\tilde{U}\subset U,|\tilde{U}|=k}}\P(\cmax=U,\Cmax=\tilde{U})\\[0.3cm]
&\leq \sum_{m\geq \frac{n}{t\log n}}Cm^{-1+{\max\{\epsilon+\delta,6\delta\}}}\sum_{m^{\frac{1}{2}-\delta}\le k\le m^{\frac{1}{2}+\delta}} \P(\cE_{m,k})\\[0.3cm]
&\leq C''\left(\frac{n}{\log n}\right)^{-1+\max\{\epsilon+\delta, 6\delta\}} \sum_{m,k}\P(\cE_{m,k})\leq C''\left(\frac{n}{\log n}\right)^{{-1+\max\{\epsilon+\delta, 6\delta\}}}.
\end{align*}
For the last bound, we have used the fact that $\max\{\epsilon +\delta ,6\delta\}<1$, as per \eqref{eq:epsdeletaprop}.\qedhere
\end{proof}

\section{Extension to General Weights}\label{sec:conc}
In our bounds for the error accumulated at each step of repair, we conditioned on a specific outcome of $\vec{K}_{n}^{1}$. In particular, for a given vertex, we conditioned on the knowledge of which one of its outgoing edges has minimal weight. When the edge weights are i.i.d.\ exponential, the memoryless property makes this conditioning particularly tractable. That being said, all our error bounds can be adapted to the setting of general weights subject to the hypothesis in Theorem \ref{thm:main}. To reiterate, the key features of the exponential family that we use are the following. Let $\{Y_{i}\}_{i=1}^{n}$ be a family of i.i.d.\ Exponential$(1)$ random variables. Then 
\begin{enumerate}
\item If $J=\text{argmin}_{1\leq i\leq n} \{Y_{i}\}$ then $Y_{J}$ is exponentially distributed with mean $1/n$.

\item Conditioned on $J$, the random variables $Y_{k}-Y_{J}$ for $k\neq J$ are i.i.d.\ Exponential$(1)$.
\end{enumerate}
We will build analogues of the above properties in the following lemmas, that extend the exponential case to the general family of random variables. Note that this family includes uniform random variables, as well as the absolute value of $\alpha-$stable random variables for $\alpha\in (1,2]$. \\

In the following lemma we study the asymptotic properties of the minimum of $m$ random variables. 
\begin{lemma}\label{lem:genweight}
Let $X_1, X_2, \cdots, X_m$ be non-negative i.i.d.\ continuous random variables with density function $p$, which has the following properties:
\begin{itemize}
    \item $p$ is continuous at $0$,
    \item $p(0)>0$,
    \item $\displaystyle\int_0^{\infty}xp(x)dx<\infty$.
\end{itemize}
Then the following hold. 
\begin{enumerate}
\item For all integers $k\ge 1$, the $k$-th moment of the minimum satisfies 
\[
\lim_{m\to\infty}m^k\E(\min\{X_1, X_2, \cdots, X_m\}^k)=\frac{k!}{p(0)^k}.
\]
\item Let $q^{(m)}$ denote the density of $m\cdot \min\{X_1, X_2, \cdots, X_m\}$. Then $q^{(m)}(0)=p(0)$, and there exist $\eps,c>0$ such that $q^{(m)}(x)\geq c$ for all $x\in [0,\eps]$. 
\end{enumerate}
\end{lemma}

\begin{proof}

Let $F(x)$ be the distribution function of $X_1$ and $Y_m = \min\{X_1, X_2, \cdots, X_m\}$. We can express the moments of $Y_m$ in terms of its tails to get
\[
  \E(Y_m^k)=\int_0^\infty kx^{k-1} \P(Y_m > x)dx = \int_0^\infty kx^{k-1} (1-F(x))^m dx.
\]
Assuming $t=mx$, we get
\begin{equation}\label{eq:lim_moment}
m^k\E(Y_m^k)=k\int_0^{\infty}t^{k-1}\left(1-F\left(\frac{t}{m}\right)\right)^m dt.
\end{equation}
Now using a Taylor expansion of $F$ up to the linear term, and since $F(0)=0$ we have
\[
F(x)=p(0)x+o(x)
\]
as $x \downarrow 0$. For any $t\ge 0$,
$$mF\left(\frac{t}{m}\right)\to p(0)t,$$
as $m\to\infty$. Consequently, as $m\to\infty$,
$$\left(1-F\left(\frac{t}{m}\right)\right)^m\to e^{-p(0)t}.$$
Thus the integrand in \ref{eq:lim_moment} converges pointwise to $t^{k-1}e^{-p(0)t}$. From the continuity of $p$ at $0$ we can consider a $\delta>0$, such that for all $0\le x \le \delta$,
\[
F(x) \ge \frac{p(0)}{2}x.
\]
Now
\begin{equation}\label{eq:cutoff}
 m^k\E(Y_m^k)= \int_0^\delta m^kkx^{k-1} \left( 1 - F(x) \right)^m dx + \int_\delta^\infty m^kkx^{k-1} (1-F(x))^m dx. 
\end{equation}

For the second term of \eqref{eq:cutoff} we have, if $m>k$,
\begin{align*}
   \int_\delta^\infty m^kkx^{k-1} (1-F(x))^m dx & \le m^k(1-F(\delta))^{m-k} \int_\delta^\infty kx^{k-1}(1-F(x))^k dx \\[0.3cm]
   & \le m^k(1-F(\delta))^{m-k} \,\E(Y_k^k)
\end{align*}
where $Y_k$ is the minimum of $k$ independent variables with density $p$.  Notice this term decays exponentially in $m$ because $F(\delta)>0$ and $\E(Y_k^k) \leq  \prod_{i=1}^k \E(X_i) = \E(X_1)^k$ which is finite by hypothesis. When $0\le t\le m\delta$, $(1-F(t/m))^m\le e^{-p(0)t/2}$ and $t^{k-1}e^{-p(0)t/2}$ is integrable on $[0,\infty]$. By dominated convergence theorem and the fact that $(1-F(t/m))^m\to e^{-p(0)t}$, we obtain
$$\lim_{m\to\infty}m^k\E(Y_m^k)=\lim_{m\to\infty}k\int_0^{m\delta}t^{k-1}\left(1-F\left(\frac{t}{m}\right)\right)^m dt=k\int_0^{\infty}t^{k-1}e^{-p(0)t}dt=\frac{k\Gamma(k)}{p(0)^k}=\frac{k!}{p(0)^k}.$$
As for the density, let $r$ denote the density of $\min\{X_1, X_2, \cdots, X_m\}$. Then 
\begin{align*}
r(t)&=-\frac{d}{dt}\P\big(\min\{X_1, X_2, \cdots, X_m\}>t\big)\\
&=-\frac{d}{dt}\big(\P(X_1>t)\big)^{m}=-m\big(\P(X_1>t)\big)^{m-1}\frac{d}{dt}\P(X_1>t)=m\big(\P(X_1>t)\big)^{m-1}p(t)
\end{align*}
Next, note that \[
q^{(m)}(t)=\frac{1}{m}r\left(\frac{t}{m}\right)=\bigg(\P\big (X_{1}>\frac{t}{m}\big)\bigg)^{m-1}p(t).
\] 
Since $p$ is continuous at $0$, there exists $\eps>0$ such that $p(t)>\frac{p(0)}{2}$ for all $t\in [0,\eps]$. Thus, for such $t$, 
\[
q^{(m)}(t)\geq \left(1-\int_{0}^{\eps/m}p(t)dt\right)^{m-1}p(t)\geq \left(1-\frac{p(0)\eps}{2m}\right)^{m-1}\frac{p(0)}{2}\geq \frac{p(0)}{4}e^{-\frac{p(0)\eps}{2}}.
\]
\end{proof}
We can also get the centred moments.
\begin{remark}
Using the same techniques as Lemma \ref{lem:genweight}, we can deduce that
$$\E(Y_m^k)\sim\frac{k!}{m^kp(0)^k}.$$
By expanding the binomial and using linearity of expectation we get
\[
   \E((Y_m-\E(Y_m))^k) \sim \frac{c_k}{m^k p(0)^k}\text{ where } c_k=k!\sum_{j=0}^k \frac{(-1)^j}{j!}.
\]

\end{remark}

\begin{lemma}\label{lem:ind_gen}
Let $X_1,X_2,\cdots,X_m$ be non-negative independent random variables with densities $p_1,\cdots,p_m$ respectively and finite expectation. Moreover, assume that $p_i$ are continuous at $0$, and that exist $\eps>0$ and $c>0$, such that 
$$p_i(x)>c\text{ for all }x\in[0,\eps] \text{ and } \E (X_{i})\text{ uniformly bounded}$$
for all $1\le i\le m$, then
$$\E(\min\{X_1,\cdots,X_m\})=O\left(\frac{1}{m}\right).$$
\end{lemma}
\begin{proof}
Let
$$M_m:=\min\{X_1,\cdots,X_m\},$$
and $F_i$ be the distribution function of $X_i$. Now
$$\E(M_m)=\int_0^{\infty}\P(M_m>x)dx=\int_0^{\infty}\prod_{i=1}^m(1-F_i(x))dx,$$
by independence.
Since there are finitely many random variables, the continuity of each $p_i$ at $0$ implies that there is some $\eps >0$, such that 
$$p_i(x)>c,$$
for all $x\in [0,\eps]$ and $1\le i\le m$. So for $0\le x\le\eps$, we have $F_i(x)\ge cx$, for each $i$ and 
$$\prod_{i=1}^m(1-F_i(x))\le (1-cx)^m\le e^{-mcx}.$$
For $x>\eps$, the monotonicity of $F_i$ and the fact that $F_i(\eps)\ge c\eps$ give us
$$\prod_{i=1}^m(1-F_i(x))\le(1-F_1(x))\prod_{i=2}^m(1-F_i(\eps))\le (1-c\eps)^{m-1}(1-F_1(x)).$$
Now
$$\int_0^{\eps}\prod_{i=1}^m(1-F_i(x))dx\le\int_{0}^{\eps}e^{-mcx}dx\le\frac{1}{mc},$$
and
$$\int_{\eps}^{\infty}\prod_{i=1}^m(1-F_i(x))dx\le (1-c\eps)^{m-1}\int_{\eps}^{\infty}(1-F_1(x))dx\le (1-c\eps)^{m-1}\E(X_1).$$
Combining these two and the fact that $0<1-c\eps<1$, we can conclude that
\[
\E(M_m)=O\left(\frac{1}{m}\right).\qedhere
\]

\end{proof}

For the following lemma we also consider the same hypotheses as in Lemma \ref{lem:genweight}.
\begin{lemma}\label{lem:genmemoryless}
    Let $J=\text{argmin}_{1\leq i\leq m} \{X_{i}\}$. Conditioned on $J$ and $X_J=x$, for almost all $x$ such that $p(x)>0$ and $F(x)<1$, the random variables $\{X_i - X_J \}_{i\neq J}$ are independent and have common density 
    \[
      p_x(t)=\frac{p(x+t)}{1-F(x)}, \quad t\geq 0,
    \]
where $p$ is the density of $X_i$ and $F$ is distribution function of $X_i$. 
\end{lemma}
\begin{proof}
    We will show that, conditional on $J=j$ and on $X_J=x$, the joint density of $\{X_i - x\}_{i\neq j}$ is the product of the one-dimensional densities.\\

    First notice that, from the independence, the joint density of all the variables is 
    \[
      f_{X_1,\cdots,X_m}(x_1,\cdots,x_m)=\prod_{i=1}^m p(x_i).
    \]
    Conditioned on the index where the minimum is attained $J=j$, and the minimum $X_J=x$ we can make the change of variables
    \[
      R_i = X_i - x \quad \text{for } i\neq j
    \]
    which has Jacobian equal to 1.\\

    The joint distribution of $(J, X_J)$ and $\{X_i\}_{i\neq J}$ can be computed using this change of variables to get
    \[
    f_{J,X_J,\{X_i\}_{i\neq J}}(j,x,(x_i)_{i\neq J}) = p(x)\prod_{i\neq j} p(x+r_i).
    \]

    The marginal distribution of $(J, X_J)$ is straightforward to compute: we integrate the variables that are not the minimum and consider our change of variables. With this we have
    \begin{align*}
        f_{J,X_J}(j,x) & =  p(x)\prod_{i\neq j} \int_0^\infty p(x+r_i) dr_i\\[0.2cm]
        & =p(x)\left( \int_x^\infty p(r) dr\right)^{m-1}\\[0.4cm]
        & =p(x)\left( 1-F(x)\right)^{m-1}.
    \end{align*}
    Finally, we consider the quotient to get the conditional density of the variables $R_i = X_i - x$:
    \[
     f_{(R_i)_{i\neq j}|X_J=x, J=j}((r_i)_{i\neq j})=\frac{p(x)\prod_{i\neq j} p(x+r_i)}{p(x)(1-F(x))^{m-1}} = \prod_{i\neq j}\frac{p(x+r_i)}{1-F(x)}. \qedhere
    \]
    
\end{proof}
The following corollary is immediate. 
\begin{corollary}\label{cor:expbound}
Let $J=\text{argmin}_{1\leq i\leq m} \{X_{i}\}$. Conditioned on $J$ and $X_J=x$ as per previous lemma for all $i\neq J$, 
\[\E[(X_i - X_J)|J, X_{J}=x]\leq \frac{\E (X_{1})}{1-F(x)}.\]    
\end{corollary}
\begin{proof}
By Lemma \ref{lem:genmemoryless}, 
\[
\E[(X_i - X_J)|J, X_{J}=x]=\int_{0}^{\infty}t\frac{p(x+t)}{1-F(x)}dt\leq \frac{1}{1-F(x)}\left(\int_{x}^{\infty}yp(y)dy-x\int_{x}^{\infty}p(y)dy\right). 
\]
The difference of the integrals is bounded above by $\E (X_{1})$. 
\end{proof}
From the previous lemma, we can derive a uniform lower bound to finally prove that the variables $\{ X_i-X_J\}_{i\neq J}$ still have expectation $O(1/m)$. This allows us to extend our repair technique to these more general weights.

\begin{lemma}\label{lem:minpdflb}
Let $p(x)$ and $p_x(t)$ be densities defined as in the previous lemmas. Then there exist constants $\eps, c >0$ such that if the minimum $X_J$ is in the interval $[0, \eps]$
\[
  p_x(t)>c
\]
for almost every $x\in [0,\eps]$ and $0\le t\le\eps$.
\end{lemma}
\begin{proof}

From the continuity of $p(x)$ at $0$ and the fact that $p(0)>0$ we can find a $\eps>0$ such that $p(y)>\frac{p(0)}{2}>0$ for every $y\in [0,2\eps]$. This implies that if the minimum $X_J=x$ is in $[0,\eps]$ then 
\[
  p_x(t)=\frac{p(x+t)}{1-F(x)} \geq p(x+t) > \frac{p(0)}{2}=c,
\]
for all $t\in [0,\eps]$.\qedhere
 
\end{proof}

With the previous Lemma, we can apply the same technique as in Lemma \ref{lem:genweight} to bound the expectation of the minimum and finally conclude. 

\begin{lemma} \label{lem:unifbound}
    Let $X_1, X_2, \cdots, X_m$ be non-negative i.i.d.\ continuous random variables with density function the same as in Lemma \ref{lem:genweight}, then
\[
\E (\min_{i\neq J} (X_i-X_J)|J) = O\left(\frac{1}{m}\right).
\]
\end{lemma}
\begin{proof}
Let
$$Z_m=\min_{i\neq J}(X_i-X_J).$$
We consider two cases separately.

\textsf{Case 1 $(x\le\eps)$}. Let $G_x$ be the common distribution function of $(X_i-X_J)$ for all $i\neq J$. Then by Lemma \ref{lem:minpdflb}, for all $t\in [0,\eps]$,
\[
G_x(t)\ge ct.
\]
Therefore,
\[
\P(Z_m>t|J,X_J=x)=(1-G_x(t))^{m-1}\le e^{-(m-1)ct},\hspace{2mm}0\le t\le\eps.
\]
For $t>\eps$, by monotonicity, we get
\[
(1-G_x(t))^{m-1}\le (1-c\eps)^{m-2}(1-G_x(t)).
\]
Note that
$$\int_0^{\infty}(1-G_x(t))dt=\frac{\int_x^{\infty}(1-F(u))du}{1-F(x)}\le\frac{\E(X_1)}{1-F(\eps)}.$$
The last inequality comes from the fact that $\int_x^{\infty}(1-F(u))du\le\int_0^{\infty}(1-F(u))du=\E(X_1)$. Now, for almost all $x\in [0,\eps]$,
\begin{multline*}
\E(Z_m|J,X_J=x)=\int_0^{\eps}\P(Z_m>t|J,X_J=x)dt+\int_{\eps}^{\infty}\P(Z_m>t|J,X_J=x)dt\\\le\int_0^{\eps}e^{-(m-1)ct}dt+ (1-c\eps)^{m-2}\frac{\E(X_1)}{1-F(\eps)}=O\left(\frac{1}{m}\right).
\end{multline*}
Integrating over $X_J\le\eps$ yields $\E(Z\mathds 1_{X_J\le\eps})=O(1/m)$.\\

\textsf{Case 2 $(x>\eps)$}. We know that $Z_m\le\sum_{i=1}^m X_i$. Thus we obtain
\[
\E\left(Z_m\mathds1_{X_J>\eps}\right)\le m\E\left(X_1\mathds1_{X_1>\eps}\right)(1-F(\eps))^{m-1}=O(1/m).
\]
Since $J$ is uniform and independent of ordered sampled values, we have
\[
\E(Z_m|J)=\E(Z_m)=O(1/m).\qedhere
\]
\end{proof}

\begin{remark} \label{rem:minbound}
    Notice that in all our repair steps, we are optimizing over many vertices and not just one. To justify the application of the previous lemma for repairs in this general-weight setting, we further require that all of the minima are simultaneously small, in the small interval $[0,\eps]$. This follows directly from a union bound. Let $A\subseteq V_n\setminus \{ \partial \}$, then
    \[
    \P\left(\max_{v\in A} X_v^{(m)}>\eps\right)\le |A|(1-F(\eps))^m
    \]
    which decays to $0$ exponentially fast.
\end{remark}

We now verify an analogous result to Lemma \ref{lem:Wclt} for the general family of weights. Combined with the previous lemmas in this section we can extend our arguments for exponential random variables to the general class of i.i.d. edge weights, and verify Theorem \ref{thm:main}. 

\begin{lemma}
    Let $W_m = \sum_{i=1}^m Y_{m,i}$ where $\{Y_{m,i} \}_{1\le i\le m}$ is a collection of independent copies of the variable $Y_m = \min\{X_1, X_2, \cdots, X_m\}$, then
    \[
    \sqrt{m} \left( \sum_{i=1}^m Y_{m,i} - m\E(Y_{m,1}) \right) \overset{d}{\rightarrow} \cN\left(0,\frac{1}{p(0)^2}\right).
    \]
\end{lemma}
\begin{proof}
To obtain a CLT for $W_m$ it is enough to check that Lyapunov's condition holds, namely:
   \[
   \lim_{m\rightarrow \infty} \frac{\sum_{i=1}^{m}\E((Y_{m,i}-\E(Y_{m,i}))^4)}{\left( \Var(\sum_{i=1}^{m}Y_{m,i}) \right)^2} = 0.
    \]
   By Lemma \ref{lem:genweight}, 
    \[
     \E((Y_{m,1}-\E(Y_{m,1}))^4)\sim\frac{9}{m^4 p(0)^4} 
    \]
    and 
    \[
    \Var(Y_{m,1}) \sim \frac{1}{m^2 p(0)^2}, 
    \]
    so therefore 
    \[
     \frac{\sum_{i=1}^{m}\E((Y_{m,i}-\E(Y_{m,i}))^4)}{\left( \Var(\sum_{i=1}^{m}Y_{m,i}) \right)^2}=\frac{\E((Y_{m,1}-\E(Y_{m,1}))^4)}{m\left( \Var(Y_{m,1}) \right)^2} \sim \frac{9}{m} \rightarrow 0.
    \]
   We can conclude that
    \[
     \frac{\sum_{i=1}^m Y_{m,i} - m\E(Y_{m,1})}{\sqrt{m\Var(Y_{m,1})}} \overset{d}{\rightarrow} \cN(0,1).
    \]
    To simplify the expression we can consider
    \[
    \sqrt{m} \left( \sum_{i=1}^m Y_{m,i} - m\E(Y_{m,1}) \right) = m\sqrt{\Var (Y_{m,1})} \cdot \frac{\sum_{i=1}^m Y_{m,i} - m\E(Y_{m,1})}{\sqrt{m\Var(Y_{m,1})}}
    \]
    and from the asymptotics of the variance, the first factor is asymptotically $\frac{1}{p(0)}$. Thus, 
    \[
    \sqrt{m} \left( \sum_{i=1}^m Y_{m,i} - m\E(Y_{m,1}) \right) \overset{d}{\rightarrow} \cN\left(0,\frac{1}{p(0)^2}\right). \qedhere
    \]
\end{proof}

We have assembled all the ingredients needed to prove Theorem \ref{thm:main} in full generality. Note, we only need to establish that Lemmas \ref{lem:errbd1}, \ref{lem:errbd2} and \ref{lem:errbd3} hold when our weights belong to the general class of distributions. In addition to the graph structure of $\vec{K}_{n}^{1}$, we include weights of the edges in $\vec{K}_{n}^{1}$ in the conditioning. We restrict to the intersection of the good event $\cE_{\epsilon,\delta,\eta}$, and $\{\max_{v\in V_{n}}X^{(n-1)}_{v}<\eps\}$ where $X_{v}^{(n-1)}$ is the weight of the outgoing edge at $v$ in $\vec{K}^{1}_{n}$ (the minimum of $n-1$ i.i.ds), and $\eps$ is as per Remark \ref{rem:minbound}. Clearly, the intersection of the two events hold with high probability. Note that each of our repair steps follow the same basic template. We begin with a collection of candidate repair sites. For a given site, situated at a vertex $u$, we consider the cheapest repair that can be executed at that site . Finally, we optimize over all candidate $u$. A crucial point, common to all three steps, is that at this first stage, identifying possible repairs at $u$, the number of repairs available is the same for all $u$! We illustrate the method for {\sf Step 1}, specifically the adaptation of Lemma \ref{lem:errbd1}.  Consider $u \in {\sf Cycle}_{i}$, and the collection of all edges from $u$ to $\cmax$, which we denoted $E_{\max}(u)$. Clearly, $|E_{\max}(u)|=|\cmax|$. On the event $\cE_{\epsilon,\delta,\eta}\cap \{\max_{v\in V_{n}}X^{(n-1)}_{v}<\eps\}$, optimal repair cost at $u$, given by $X_{u}:=\min \{Y_{\vec{e}}-Y_{\vec{e}_{*}(u)}:\vec{e}\in E_{\max}(u)\}$ satisfies the hypotheses of Lemmas  \ref{lem:genmemoryless} and \ref{lem:minpdflb}. Finally, if we consider $\min\{|\cmax|\cdot X_{u}:u\in {\sf Cycle}_{i}\}$, this family satisfies the hypotheses of   Lemma \ref{lem:ind_gen}, in particular using point 2 of Lemma \ref{lem:genweight}. Thus, the expectation bound of Lemma \ref{lem:errbd1} can be recovered. The same points are applicable to {\sf Step 2} and {\sf Step 3}. 

\section{Concluding Remarks}\label{sec:concreal}
We take a moment to explain  how Theorem \ref{thm:main} and the fact $\sqrt n(W_n^*-W_n)\to 0$ in probability help us obtain the CLT of \cite{hansen1997limit} for the  family of weights considered in this article. Let $T_n$ be the minimum spanning arborescence on $\Kn$ \textit{without a} fixed boundary and let $c(T_n)$ denote its cost. Then the CLT of \cite{hansen1997limit} (paraphrased) states that
$$\sqrt n(c(T_n)-n\E(X_v))\xrightarrow[]{d}\cN(0,1/p(0)^2).$$

This can be extracted from the results of this paper, provided we are able to show that
$$\sqrt n(W_n^*-c(T_n))\to 0,$$
in probability. Let $\rho$ be the boundary of $T_n$. Since $c(T_n)$ minimizes the cost among all spanning arborescences (with all possible boundaries), we have
$$W_n^*\ge c(T_n).$$
Recall that $X_v$ is the weight of minimal outgoing edge from $v\in V_n$. Note, $c(T_n)\ge\sum_{v\neq\rho}X_v$ (since every vertex apart from $\rho$ must have exactly one out going edge),
$$
W_n^*-c(T_n)\le W_n^*-\sum_{v\neq\rho}X_v=W_n^*-(W_n-X_{\rho}+X_{\partial})\le W_n^*-W_n+\max_{v\in V_n}X_v.
$$
Next, we demonstrate that $\sqrt n\max_{v\in V_n}X_v\to 0$ in probability. Since each $X_v$ is the minimum of $(n-1)$ random variables with finite expectation, 
\[
\E (X_{v}) \sim cn^{-1}\text{ and }\E(X_{v}-\E (X_{v}))^{4}\leq Cn^{-4}.
\]
Indeed, these follow from Lemma \ref{lem:genweight}. Thus, for any $\eps>0$
\[
\P(X_{v}>\eps n^{-\frac{1}{2}})\leq \P(|X_{v}-\E (X_{v})|>\eps {n}^{-\frac{1}{2}}-cn^{-1})\leq \P(|X_{v}-\E (X_{v})|>\frac{1}{2}\eps {n}^{-\frac{1}{2}}),
\]
final inequality holding for $n$ large enough.  Applying a union bound, we obtain
$$\P\left(\sqrt{n}\max_{v\in V_n}X_v>\eps\right)\le\sum_{v\in V_n}\P(X_v>\eps n^{-\frac{1}{2}})\le n \cdot \frac{2^{4}\cdot Cn^{-4}}{\eps^{4} \cdot n^{-2}}\leq \frac{2^{4}\cdot C}{\eps^{4}}\cdot n^{-1}.$$
Therefore $\sqrt n\max_{v\in V_n}X_v\to 0$ in probability. Combining the facts $0\le W_n^*-c(T_n)\le W_n^*-W_n+\max_{v\in V_n}X_v$ and $\sqrt n(W_n^*-W_n)\to 0$ in probability, we can conclude that $\sqrt n(W_n^*-c(T_n))\to 0$ in probability. 
\\

We expect that our method generalizes to families of non-negative random variables other than those with densities continuous at $0$, satisfying $p(0)\in(0,\infty)$. In particular, we anticipate that only the scaling changes, based on the application of Laplace asymptotics. If $p(x)\sim x^{\beta}$ near zero for $\beta\in (-1,\infty)$, mirroring the sequence of steps in Lemma \ref{lem:genweight}, the minimal weight edge emanating from each vertex is now $O\left(n^{-\frac{1}{1+\beta}}\right)$. The variance also transforms to $O\left(n^{-\frac{2}{1+\beta}}\right)$. The Lyapunov condition also holds, yielding a CLT for $W_{n}$ with the new scaling. However, subtleties could arise in verifying analogues of Lemmas \ref{lem:genmemoryless},  \ref{lem:minpdflb}, and \ref{lem:unifbound} . Where our methods possibly fail completely are densities which vanish to infinite order at $0$, an example given by the canonical $C^{\infty}$ bump function on $[0,1]$. Our methods are quite robust with respect to tail conditions on $Y_{\vec{e}}$, power decay for $\P(Y_{\vec{e}}>t)$ is completely admissible, even without finite first moment. However, all our estimates break with pathological tail examples such as $\P(Y_{\vec{e}}>t)\sim \frac{C}{\log t}$ as $t\to \infty$. 

\nocite{*}
\addcontentsline{toc}{section}{References}
\printbibliography

\end{document}